\documentclass[a4paper,11pt,twoside,reqno]{amsart}
\usepackage{amscd,amssymb,amsxtra}
\usepackage[mathcal]{euscript}		   

\makeatletter

\def\@maketitle{%
  \normalfont\normalsize
  \let\@makefnmark\relax  \let\@thefnmark\relax
  \ifx\@empty\@date\else \@footnotetext{\@setdate}\fi
  \ifx\@empty\@subjclass\else \@footnotetext{\@setsubjclass}\fi
  \ifx\@empty\@keywords\else \@footnotetext{\@setkeywords}\fi
  \ifx\@empty\thankses\else \@footnotetext{%
    \def\par{\let\par\@par}\@setthanks}\fi
  \@mkboth{\@nx\shortauthors}{\@nx\shorttitle}%
  \global\topskip42\p@\relax 
  \@settitle
  \ifx\@empty\authors \else \@setauthors \fi
  \ifx\@empty\@dedicatory
  \else
    \baselineskip18\p@
    \vtop{\centering{\footnotesize\itshape\@dedicatory\@@par}%
      \global\dimen@i\prevdepth}\prevdepth\dimen@i
  \fi
\box\abbox
  \@setabstract
  \normalsize
  \if@titlepage
    \newpage
  \else
    \dimen@34\p@ \advance\dimen@-\baselineskip
    \vskip\dimen@\relax
  \fi
} 
\def\enddoc@text{
}
\def\@sect#1#2#3#4#5#6[#7]#8{%
  \edef\@toclevel{\ifnum#2=\@m 0\else\number#2\fi}%
  \ifnum #2>\c@secnumdepth \let\@secnumber\@empty
  \else \@xp\let\@xp\@secnumber\csname the#1\endcsname\fi
 \ifnum #2>\c@secnumdepth
   \let\@svsec\@empty
 \else
    \refstepcounter{#1}%
    \edef\@svsec{\ifnum#2<\@m
       \@ifundefined{#1name}{}{%
         \ignorespaces\csname #1name\endcsname\space}\fi
       \@nx\textup{%
      \@nx\bfseries 
         \csname the#1\endcsname.}\enspace
    }%
  \fi
  \@tempskipa #5\relax
  \ifdim \@tempskipa>\z@ 
    \begingroup #6\relax
    \@hangfrom{\hskip #3\relax\@svsec}{\interlinepenalty\@M #8\par}%
    \endgroup
    \ifnum#2>\@m \else \@tocwrite{#1}{#8}\fi
  \else
  \def\@svsechd{#6\hskip #3\@svsec
    \@ifnotempty{#8}{\ignorespaces#8\unskip
       \@addpunct.}%
    \ifnum#2>\@m \else \@tocwrite{#1}{#8}\fi
  }%
  \fi
  \global\@nobreaktrue
  \@xsect{#5}}
\def\section{\@startsection{section}{1}%
  \z@{.7\linespacing\@plus\linespacing}{.5\linespacing}%
  {\normalfont\fontsize{12}{14}\selectfont\bfseries\centering}}
\def\subsection{\@startsection{subsection}{2}%
  \z@{.5\linespacing\@plus.7\linespacing}{1.5ex \@plus .2ex}%
  {\normalfont\bfseries}}
\def\subsubsection{\@startsection{subsubsection}{3}%
  \z@{.5\linespacing\@plus.7\linespacing}{-.5em}%
  {\normalfont\bfseries}}
\def\ps@firstpage{\ps@plain
  \def\@oddfoot{\normalfont\hfill\nine[\thepage]\hfill
     \global\topskip\normaltopskip}%
  \let\@evenfoot\@oddfoot
  \def\@oddhead{\phead\hss}%
  \let\@evenhead\@oddhead}
\newtheoremstyle{boldhead}
{\thm@preskip}
{\thm@postskip}
{\slshape}
{}
{\bfseries}
{.}
{ }
{\thmname{#1}\thmnumber{ #2}\thmnote{ (#3)}}
\newtheoremstyle{boldremark}
{\thm@preskip}
{\thm@postskip}
{\upshape}
{}
{\bfseries}
{.}
{ }
{\thmname{#1}\thmnumber{ #2}\thmnote{ (#3)}}
\renewenvironment{thebibliography}[1]{%
  \@xp\section\@xp*\@xp{\refname}%
  \normalfont\nine\labelsep .5em\relax
  \renewcommand\theenumiv{\arabic{enumiv}}\let\p@enumiv\@empty
  \list{\@biblabel{\theenumiv}}{\settowidth\labelwidth{\@biblabel{#1}}%
    \leftmargin\labelwidth \advance\leftmargin\labelsep
    \usecounter{enumiv}}%
  \sloppy \clubpenalty\@M \widowpenalty\clubpenalty
  \sfcode`\.=\@m
}{%
  \def\@noitemerr{\@latex@warning{Empty `thebibliography' environment}}%
  \endlist
}
\font\sevenrm=cmr7

\newbox\fpbox
\newbox\adbox
\newbox\abbox
\setbox\abbox=\vbox{\vglue18pt}
\def\address#1{\setbox\adbox=\hbox{\let\\=\cr
\nine\baselineskip12pt\vbox{\itshape\tabskip 0pt plus15cc
\halign to\hsize{\hfil\ignorespaces {##}\hfil\cr#1\cr}}}%
\global\setbox\abbox=\vbox{\unvbox\abbox\box\adbox\vskip16pt}}

\makeatother

\DeclareMathOperator{\Bicoalg}{Bicoalg}

\DeclareMathOperator{\Coalgsq}{Coalgsq}

\DeclareMathOperator{\coev}{coev}
\DeclareMathOperator{\coim}{coim}
\DeclareMathOperator{\Coim}{Coim}
\DeclareMathOperator{\coker}{coker}
\DeclareMathOperator{\Coker}{Coker}

\DeclareMathOperator{\comod}{-comod}
\DeclareMathOperator{\Comod}{-Comod}
\DeclareMathOperator{\End}{End}
\DeclareMathOperator{\ev}{ev}
\DeclareMathOperator{\Hom}{Hom}
\DeclareMathOperator{\id}{id}
\DeclareMathOperator{\Id}{Id}

\DeclareMathOperator{\im}{im}
\DeclareMathOperator{\IIm}{Im}
\renewcommand{\Im}{\IIm}
\DeclareMathOperator{\ind}{ind}
\DeclareMathOperator{\Ind}{Ind}
\DeclareMathOperator{\kernel}{ker}
\DeclareMathOperator{\Ker}{Ker}

\DeclareMathOperator{\modul}{-mod}
\DeclareMathOperator{\modur}{mod-}
\DeclareMathOperator{\Mor}{Mor}
\DeclareMathOperator{\Ob}{Ob}
\DeclareMathOperator{\opp}{op}
\newcommand{\op}{{\opp}}

\DeclareMathOperator{\vect}{-vect}
\DeclareMathOperator{\Vect}{-Vect}

\let\<=\langle
\let\>=\rangle
\let\e=\varepsilon
\let\ge=\geqslant
\let\kk=\Bbbk
\let\le=\leqslant
\let\om=\omega
\let\ot=\circledast

\let\tens=\otimes
\let\und=\underline
\let\wb=\overline
\let\wh=\widehat
\let\xla=\xleftarrow
\let\xra=\xrightarrow

\newcommand{\rEq}{=}

\newcommand{\ca}{{\mathcal A}}
\newcommand{\cb}{{\mathcal B}}
\newcommand{\cc}{{\mathcal C}}

\newcommand{\co}{{\mathcal O}}
\newcommand{\cp}{{\mathcal P}}

\newcommand{\cu}{{\mathcal U}}
\newcommand{\cv}{{\mathcal V}}

\newcommand{\fa}{{\mathfrak A}}
\newcommand{\fm}{{\mathfrak M}}

\newcommand{\ahat}{{\widehat{\mathcal A}}}

\newcommand{\bhat}{{\widehat{\mathcal B}}}
\newcommand{\chat}{{\widehat{\mathcal C}}}

\newcommand{\Ad}[1]{$A$\nobreakdash-\hspace{0pt}}
\newcommand{\Bd}[1]{$B$\nobreakdash-\hspace{0pt}}
\newcommand{\Cbar}[1]{$\bar C$\nobreakdash-\hspace{0pt}}
\newcommand{\Cd}[1]{$C$\nobreakdash-\hspace{0pt}}
\newcommand{\Hd}[1]{$H$\nobreakdash-\hspace{0pt}}
\newcommand{\ii}{\ddot\imath}
\newcommand{\kd}[1]{$\Bbbk$\nobreakdash-\hspace{0pt}}

\newcommand{\Rd}[1]{$R$\nobreakdash-\hspace{0pt}}
\newcommand{\SSS}{{\mathfrak S}}

\newcommand{\vto}[1]{{{\mathcal V}^{\tens #1}}}
\newcommand{\vhat}{{\widehat{\mathcal V}}}
\newcommand{\vtahat}{{\widehat{\mathcal V\tens\mathcal A}}}
\newcommand{\vtvhat}{{\widehat{\mathcal V\tens\mathcal V}}}
\newcommand{\vta}{{\mathcal V\tens\mathcal A}}
\newcommand{\vtv}{{\mathcal V\tens\mathcal V}}

\newcommand{\barten}{\bar\otimes}
\newcommand{\colimit}[1]{\lim_{\substack{\longrightarrow\\#1}}}
\newcommand{\limit}[1]{\lim_{\substack{\longleftarrow\\#1}}}
\newcommand{\lpti}{{}^\vee\!}
\newcommand{\pti}{^\vee}
\newcommand{\unten}{\,\underline{\!\otimes\!}\,}

\newcommand{\one}{_{(1)}}
\newcommand{\two}{_{(2)}}
\newcommand{\tre}{_{(3)}}

\newcommand{\corref}[1]{Corollary~\ref{#1}}
\newcommand{\defref}[1]{Definition~\ref{#1}}

\newcommand{\propref}[1]{Proposition~\ref{#1}}
\newcommand{\remref}[1]{Remark~\ref{#1}}
\newcommand{\secref}[1]{Section~\ref{#1}}
\newcommand{\thmref}[1]{Theorem~\ref{#1}}

\providecommand{\nine}{\fontsize{9}{11}\selectfont}

\swapnumbers

\theoremstyle{boldhead}

\newtheorem{theorem}[subsubsection]{Theorem}
\newtheorem{corollary}[subsubsection]{Corollary}

\newtheorem{proposition}[subsubsection]{Proposition}

\theoremstyle{boldremark}

\newtheorem{definition}[subsubsection]{Definition}
\newtheorem{remark}[subsubsection]{Remark}
\newtheorem{example}[subsubsection]{Example}
\newtheorem*{examples}{Examples}

\numberwithin{equation}{subsection}

\def\phead{\setbox\fpbox=\hbox{\sevenrm
************************************************}%
\noindent\vbox{\sevenrm\baselineskip9pt\hsize\wd\fpbox%
\centerline{***********************************************}

\centerline{BANACH CENTER PUBLICATIONS, VOLUME **}

\centerline{INSTITUTE OF MATHEMATICS}

\centerline{POLISH ACADEMY OF SCIENCES}

\centerline{WARSZAWA 19**}}\hfill}

\begin{document}

\title[Squared Hopf Algebras]%
{Squared Hopf algebras and reconstruction theorems}
\author[V.~V. Lyubashenko]{Volodymyr V. Lyubashenko}
\address{Mathematical Methods of Systems Analysis\\
Department of Applied Mathematics\\
Kiev Polytechnic Institute\\
37 prosp. Peremogy,
252056 Kiev,
Ukraine}

\subjclass{Primary 16W30, 18D35, Secondary 17B37, 18D10}

\thanks{q-alg/9605035}
\thanks{This paper will appear in the Proceedings of Banach Center
Minisemester on Quantum Groups, November 1995}

\thanks{The full version of this paper containing proofs
will be published elsewhere}
\thanks{Research was supported in part
by the EPSRC research grant GR/G 42976}

\begin{abstract}
Given an abelian $\kk$-linear rigid monoidal category $\cv$, where
$\kk$ is a perfect field, we define \emph{squared coalgebras} as
objects of cocompleted $\vtv$ (Deligne's tensor product
of categories) equipped with the appropriate notion of
comultiplication. Based on this, (squared) bialgebras and Hopf
algebras are defined without use of braiding. If $\cv$ is
the category of $\kk$-vector spaces, squared (co)algebras
coincide with conventional ones. If $\cv$ is braided,
a braided Hopf algebra can be obtained from a squared one.

Reconstruction theorems give equivalence of squared co- (bi-,
Hopf) algebras in $\cv$ and corresponding fibre functors to $\cv$
(which is not the case with the usual definitions).  Finally, squared
quasitriangular Hopf coalgebra is a solution to the problem of defining
quantum groups in braided categories.
\end{abstract}

\maketitle


\allowdisplaybreaks[1]

\section*{Introduction}
Classical reconstruction theorem (e.g. Saavedra
\cite[Section~2.3.2.1]{SaaRiv})
asserts that a \kd-coalgebra can be reconstructed from the underlying
functor from its category of comodules to vector spaces. Saavedra
\cite[Section~2.6.3~a)]{SaaRiv} and later Schauenburg~\cite{Sch:hopf1} also
prove that an essentially small abelian \kd-linear category equipped with
an exact faithful functor $\om$ to the category of finite dimensional
\kd-vector spaces is equivalent to the category of finite dimensional
comodules  over some \kd-coalgebra. A direct attempt to generalise these
results replacing the category of vector spaces by an abelian \kd-linear
rigid monoidal category $\cv$ fails. For instance, the category of
comodules over the coalgebra constructed from $\om$ is bigger than the
initial category (precise results are formulated by
Pareigis~\cite{Par:coend}).  However, if one modifies the definitions of
coalgebras and comodules in a monoidal category, the reconstruction theorem
will be recovered.  This is the main conclusion of this work.

The new notion will be called a {\em squared coalgebra}. The monoidal
version of the reconstruction theorem dictates the definition of a squared
bialgebra. Squared Hopf (co)algebras based on $\cv$ can be also defined,
even if $\cv$ is not braided, but satisfies a much weaker condition! If
$\cv$ is braided, a squared Hopf (co)algebra determines a braided Hopf
algebra, but not vice versa. Finally, squared quasitriangular Hopf
coalgebra is a solution to the problem of defining quantum groups in
braided categories.

All definitions are based on the notion of tensor product of abelian
categories given by Deligne \cite{Del:GF}. The squared notions (coalgebras,
bialgebras, Hopf algebras) are objects of the cocompleted tensor square of
the initial category $\cv$, whence the terminology. The structure maps --
comultiplication, multiplication etc. -- are morphisms in tensor powers of
$\cv$. The associativity and other properties mean equality of two
composite morphisms in tensor powers of $\cv$.

More precisely, we use the cocompletions of tensor powers of $\cv$, where
the cocompletion of an abelian \kd-linear category $\ca$ with finite
dimensional \kd-spaces of morphisms always means the category $\hat\ca =
\ind-\ca$, made of filtered inductive limits of objects of $\ca$. If $\ca$
is essentially small, $\ahat$ is equivalent to the category
$\und{\Hom}_{\kk,\text{l.e.}} (\ca^\op,\kk\vect)$ of left exact functors
$\ca^\op \to \kk\vect$.  (See Grothendieck and Verdier \cite{GroVer}.)

We assume that $\kk$ is a perfect field.

Let us recall the reconstruction theorem in details. If
$\omega:\cc\to \kk\vect$ is a faithful exact \kd-linear functor and
$\cc$ is essentially small, then there is an equivalence $F$ of $\cc$
with the category of \Cd-comodules, where $C$ is the \kd-coalgebra
\[ C = \int^{X\in\cc} \omega(X) \tens_\kk \omega(X)^* ,\]
and $\omega$ is isomorphic to the composite of $F$ and the underlying
functor $\cu:C\comod \to \kk\vect$. When $\kk\vect$ is replaced by an
abelian \kd-linear rigid monoidal category $\cv$ with finite
dimensional spaces $\Hom_\cv(\text-,\text-)$ such that $\End I=\kk$
($I$ is the unit object) and each object has finite length, a version of
reconstruction theorem holds, although $F$ is no longer an equivalence.
It turns out that by modifying the definitions of coalgebras and
comodules one can make $F$ into an equivalence,  thus recovering the
original form of the theorem. Namely, instead of the coalgebra in $\vhat$
\begin{equation}\label{e4a}
\bar C = \int^{X\in\cc} \omega(X) \tens \omega(X)\pti , \tag{1}
\end{equation}
one can use the squared coalgebra
\begin{equation}\label{e4b}
C = \int^{X\in\cc} \omega(X) \odot \omega(X)\pti \in\vtvhat , \tag{2}
\end{equation}
where $\odot: \cv\times\cv \to \cv\tens\cv$ is the canonical functor of
exterior tensor product.

To explain the definition of a squared coalgebra, we recall that the tensor
product functor $\tens:\cv\times\cv \to \cv$ can be decomposed as
$\cv\times\cv \xra\odot \cv\tens\cv \xra\ot \cv$ up to an isomorphism.
Using the functor $\odot$ and the diagonal restriction functor $\ot$ we can
construct various objects like $C_{13}\odot I_2 \in \wh{\cv^{\tens3}}$ ($I$
is the unit object and the subindices indicate the tensorands to which an
object is placed), or like $C_{12'}\tens C_{2''3} \in \wh{\cv^{\tens3}}$
(this is the result of applying $\ot$ on the second and third places to
$C_{12}\odot C_{34} \in \wh{\cv^{\tens4}}$), or like $C_{1'1''}= \ot C \in
\vhat$ (the dash and the double dash indicate the order of multiplicands)
etc. Notice that $\ot$ applied to (2) 
gives (1).

The reader is advised to consider the example $\cv = H\comod$ (finite
dimensional \Hd-comodules) throughout this work, where $H$ is a Hopf
\kd-algebra. Then $\vhat = H\Comod$ (all \Hd-comodules), $\cv^{\tens n} =
H^{\tens n}\comod$, $\wh{\cv^{\tens n}} = H^{\tens n}\Comod$. The functor
$\odot: \cv\times\cv \to \cv\tens\cv$ sends a pair of \Hd-comodules $(X,Y)$
to the $H\tens_\kk H$-comodule $X\tens_\kk Y$ (exterior tensor product of
comodules) and the diagonal restriction functor $\ot:\cv\tens\cv \to \cv$
sends an $H\tens_\kk H$-comodule to an \Hd-comodule via the homomorphism of
coalgebras $m: H\tens_\kk H \to H$ -- the multiplication. Thus, given an
 $H\tens_\kk H$-comodule $C$ with the coaction $\delta:  c \mapsto c\one
\tens c\two \tens c\tre \in H\tens_\kk H \tens_\kk C$ we get the
$H^{\tens3}$-comodule $C_{13}\odot I_2 = (C, \delta: c \mapsto c\one
\tens1\tens c\two \tens c\tre \in H^{\tens3}\tens_\kk C)$, the
$H^{\tens3}$-comodule $C_{12'}\tens C_{2''3} = (C\tens_\kk C, c\tens d
\mapsto c\one\tens c\two d\one\tens d\two\tens c\tre\tens d\tre \in
H^{\tens3}\tens_\kk C^{\tens2})$, the \Hd-comodule $\ot C = (C, \delta: c
\mapsto  c\one c\two\tens c\tre \in H\tens_\kk C)$ etc.

Now we may define a squared coalgebra as an object $C\in \vtvhat$ equipped
with the comultiplication $\Delta_{123}: C_{13}\odot I_2 \to C_{12'}\tens
C_{2''3}$ and the counit $\e: C_{1'1''} \to I_1$, which satisfy the
coassociativity axiom (an equation in $\wh{\cv^{\tens4}}$, see
\eqref{d23a}) and two axioms for the counit (equations
\eqref{e23b}-\eqref{e23c} in $\wh{\cv^{\tens2}}$). A \Cd-comodule is an
object $X\in\cv$ equipped with the coaction $\delta:  X_1\odot I_2 \to
C_{12'}\tens X_{2''} \in \wh{\cv^{\tens2}}$, which is coassociative
(equation \eqref{d28a} in $\wh{\cv^{\tens3}}$) and counital (equation
\eqref{e28b} in $\wh{\cv^{\tens2}}$). It turns out that for any object
$Y\in\cv$ the object $Y\odot Y\pti \in \cv\tens\cv$ has the canonical
structure of a squared coalgebra and $Y$ is a comodule over it. We deduce
that the coend (2) 
is a squared coalgebra as well.  Moreover, the second part of the
reconstruction theorem claims that any squared coalgebra is isomorphic to a
coalgebra of the form (2).

Thus, the full form of the reconstruction theorem asserts equivalence of
the following two categories: the category of \kd-linear exact faithful
functors from an essentially small category to $\cv$ and the category of
squared coalgebras in $\cv$. Philosophically, categories over the category
$\cv$ are fully encoded in terms of coalgebras living in $\cv$ (in fact, in
$\vtvhat$) and vice versa. Comparing the category of comodules ${}^C\cv$
over a squared coalgebra $C$ and the category of comodules ${}^{\bar C}\cv$
over the ordinary coalgebra (comonoid) $\bar C = \ot C$ we get ${}^{\bar
C}\cv = {}^C\cv \tens\cv$. That is expected from the description of the
coend reconstructed from the functor ${}^{\bar C}\cv \to \cv$ given by
Pareigis~\cite{Par:coend} in the case $\cv= H\comod$ for a Hopf algebra
$H$.

The monoidal version of the reconstruction theorem also holds. Namely, the
category of monoidal $\kk$ -linear exact faithful functors $\omega:
\cc\to\cv$ ($\cc$ is essentially small) and the category of squared
bicoalgebras in $\cv$ are equivalent. A {\em squared bicoalgebra} is
defined as an object of $\vtvhat$ having the structure of a squared
coalgebra and of an algebra in the monoidal category $\vtvhat$ with
compatibility axioms which require that the multiplication and the unit
were homomorphisms of coalgebras.  (There are several monoidal structures
in $\vtvhat$ and we chose a special one.)

The dual notion, squared bialgebras, is defined as a squared algebra
structure plus a coalgebra structure in $\vtvhat$ with compatibility
axioms.  Unlike the case of vector spaces the notion of squared bicoalgebra
is not self-dual, so it has to be distinguished from squared bialgebras.
The choice of terminology is motivated by our primary interest in comodules
rather than in modules. We shall simplify it further dropping the adjective
squared and keeping the term bicoalgebra.

Notice that a braiding in $\cv$ is not required for work with bicoalgebras.
However, if $\cv$ is braided any bicoalgebra $B$ generates a braided
bialgebra $\ot B$ in $\cv$. Not every braided bialgebra comes from a
bicoalgebra in that way.

To introduce Hopf algebras we require that the second dual $X^{\vee\vee}$
of an object $X\in\cv$ were isomorphic to $X$ via a functorial isomorphism
$\zeta=\zeta_X: X \to X^{\vee\vee}$. Obviously, this condition is weaker
than existence of braiding. Given a squared coalgebra $C$ in such $\cv$,
one can define the opposite coalgebra $C_\op$. A (squared) Hopf coalgebra
is defined as a bicoalgebra $H$ together with an isomorphism $\gamma: H_\op
\to H \in \vtvhat$ -- the antipode -- satisfying two equations in
$\vtvhat$. The reason for introducing Hopf coalgebras is the following: the
category of comodules over a Hopf coalgebra is rigid and the rigid version
of the reconstruction theorem holds: the category of monoidal \kd-linear
exact faithful functors $\omega: \cc \to \cv$, where $\cc$ is rigid
monoidal (and essentially small), and the category of Hopf coalgebras in
$\cv$ are equivalent. The dual notion, squared Hopf algebras, is not
equivalent to the notion of Hopf coalgebras.

If, in addition, $\cv$ is braided, we get the equivalence of the category
of monoidal \kd-linear exact faithful functors $\omega: \cc \to \cv$, where
$\cc$ is  rigid braided, and of the category of {\em quasitriangular Hopf
coalgebras}, appropriately defined. In particular, the category of
comodules over a quasitriangular Hopf coalgebra is braided. This is not
trivial, and allows to introduce a non-obvious braiding for the bigger
category of comodules over the braided Hopf algebra $\ot H$. However, it
seems impossible, in general, to introduce a braided structure of any kind
for the whole category of comodules over a braided Hopf algebra not related
with Hopf coalgebras. Thus the notion of a quasitriangular Hopf coalgebra
is the closest to the idea of a ``quantum group in a braided category''.

In particular, applying the (re)construction theorem to the identity
functor $\Id:\cv \to \cv$, we get a quasitriangular Hopf coalgebra
structure of the coend
\[ C = \int^{X\in\cv} X\odot X\pti ,\]
and this is the most interesting case for us. Similar notions exist for
ribbon categories.

\section{Tools}
\subsection{Tensor product of abelian categories}
In this work $\kk$ will denote a perfect field.

\begin{definition}
We say that an abelian \kd-linear category $\ca$ is a {\em category with
length} if
\begin{enumerate}
\item for any pair of objects $X,Y\in\ca$ the \kd-vector space
$\Hom_\ca(X,Y)$ is finite dimensional, and
\item any object $X\in\ca$ has finite length.
\end{enumerate}
\end{definition}

The most impressive result from the theory of such categories is:

\begin{proposition}[Deligne \cite{Del:GF} Corollary 2.17]\label{P10}
Let $\ca$ be an abelian \kd-linear category with length.
Assume that there exists an object $X\in\ca$ such that any object $Y\in\ca$
is a subquotient of $X^n = X\oplus\dots\oplus X$ for some $n$. Then the
category $\ca$ is equivalent to the category $\modur A$ for some finite
dimensional \kd-algebra $A$.
\end{proposition}

From now on $\kk$ is a {\em perfect} field. We shall remind the definition
of a tensor product of abelian \kd-linear categories, belonging to
Deligne \cite{Del:GF}, in a modified form, using his results, valid in
assumption of perfectness.

\begin{definition}[following Deligne \cite{Del:GF} Definition 5.1]
\label{D11}
Let $\ca_1,\dots,\ca_n$ be \kd-linear abelian categories with
 length. A \kd-linear abelian category $\ca_1\tens\dots\tens\ca_n$,
equipped with a \kd-multilinear, exact in each variable functor
\[ \odot:  \ca_1\times\dots\times\ca_n \to \ca_1\tens\dots\tens\ca_n \]
is called a tensor product of $\ca_1,\dots,\ca_n$ if for each \kd-linear
abelian category $\ca$ the induced functor
\[ \und\Hom_{\kk,\text{r.e.}} (\ca_1\tens\dots\tens\ca_n, \ca) \to
\und\Hom_{\kk,\text{r.e.}} (\ca_1\times\dots\times\ca_n, \ca), F \mapsto
F\circ\odot \]
from the category of \kd-linear right exact functors to the category of
\kd-multilinear right exact in each variable  functors is an equivalence.
\end{definition}

\begin{remark}
Equivalently, one can use left exact functors in the above definition
(Deligne \cite[Proposition 5.13]{Del:GF}).
\end{remark}

\begin{proposition}[Deligne \cite{Del:GF} Proposition 5.13]\label{P11A}
The tensor product of \kd-linear abelian categories with
length exists and is a category with length. It is unique up to an
equivalence. The functor similar to that of \defref{D11}
\[ \und\Hom_{\kk,\text{e.}} (\ca_1\tens\dots\tens\ca_n, \ca) \to
\und\Hom_{\kk,\text{e.}} (\ca_1\times\dots\times\ca_n, \ca), F \mapsto
F\circ\odot \]
with the right exact functors replaced by exact functors is also an
equivalence.  The natural map
\[ \tens_i \Hom(X_i,Y_i) \to \Hom( \odot_i X_i,\odot_i Y_i) \]
is an isomorphism.
\end{proposition}

This follows mainly from:

\begin{proposition}[Deligne \cite{Del:GF} Proposition 5.3, Corollary 5.4]
\label{P11B}
Let $\ca_i$ be $A_i\modul$, where $A_i$ are finite dimensional
\kd-algebras, $1\le i\le n$. Then $\ca = A_1\tens_\kk \dots \tens_\kk
A_n\modul$ equipped with the exterior tensor product functor $\odot:
(X_1,\dots,X_n) \mapsto X_1\tens_\kk \dots \tens_\kk X_n$ is a tensor
product of $\ca_1,\dots,\ca_n$.
\end{proposition}

Notice that any category with length $\ca$ is a filtered inductive limit of
its full subcategories of the form $\<X\>$ -- the full subcategory formed
by the objects $Y\in\ca$, which are subquotients of $X^n$ for some $n$.
Here $X$ is an object of $\ca$. This remark is used together with
Propositions \ref{P10}, \ref{P11B}.

Let $\ca_1,\dots,\ca_n,\cb_1,\dots,\cb_n$ be \kd-linear abelian categories
with length and let $T_i: \ca_i \to \cb_i$ be \kd-linear left (resp. right)
exact functors.  By definition, there exists a \kd-linear left (resp.
right) exact functor
\[ T=\tens_i T_i:  \tens_i\ca_i \to \tens_i\cb_i, \qquad T(\odot_i X_i) =
\odot_i T(X_i). \]

\begin{proposition}[Deligne \cite{Del:GF} Proposition 5.14]\label{P12.5B}
If $T_i$ are exact (resp. exact and faithful, resp. equivalence of $\ca_i$
with a full subcategory of $\cb_i$ stable with respect to subquotients),
then $T$ has the same property.
\end{proposition}

\subsection{Monoidal categories}
We shall remind only few definitions here. For the developed introduction
see \cite{Mac:work}.

\begin{definition}
A {\em rigid category} $\cc$ is a monoidal category in which for every
object $X\in\cc$ there exist dual objects $X\pti$, $\lpti X \in\cc$ and
morphisms of evaluation and coevaluation
\begin{alignat*}3
\ev &: X\tens X\pti \to I\  &=
\makebox[20mm][l]{
\raisebox{-3mm}[5mm][5mm]{
\unitlength 0.800mm
\linethickness{0.4pt}
\begin{picture}(15.94,14)
\put(10,7){\oval(10,14)[b]}
\put(4,5){\makebox(0,0)[rb]{$X$}}
\put(15.94,5){\makebox(0,0)[lb]{$X\pti$}}
\end{picture}
}} \ ,
 &\qquad \ev &: \lpti X\tens X \to I\ &=
\makebox[20mm][l]{
\raisebox{-3mm}[5mm][5mm]{
\unitlength 0.800mm
\linethickness{0.4pt}
\begin{picture}(15.94,14)
\put(10,7){\oval(10,14)[b]}
\put(4,5){\makebox(0,0)[rb]{$\lpti X$}}
\put(15.94,5){\makebox(0,0)[lb]{$X$}}
\end{picture}
}} \ ,   \\
\coev &: I\to X\pti\tens X\ &=
\makebox[20mm][l]{
\raisebox{-3mm}[5mm][5mm]{
\unitlength 0.80mm
\linethickness{0.4pt}
\begin{picture}(15.94,7)
\put(9.94,0){\oval(10,14)[t]}
\put(3.94,0){\makebox(0,0)[rb]{$X\pti$}}
\put(15.94,0){\makebox(0,0)[lb]{$X$}}
\end{picture}
}} \ ,
 &\qquad \coev &: I\to X\tens\lpti X\ &=
\makebox[20mm][l]{
\raisebox{-3mm}[5mm][5mm]{
\unitlength 0.80mm
\linethickness{0.4pt}
\begin{picture}(15.94,7)
\put(9.94,0){\oval(10,14)[t]}
\put(3.94,0){\makebox(0,0)[rb]{$X$}}
\put(15.94,0){\makebox(0,0)[lb]{$\lpti X$}}
\end{picture}
}} \ ,
\end{alignat*}
satisfying standard equations.
\end{definition}

In a rigid monoidal category $\cc$ there is a pairing
\begin{multline*}
(X\tens Y)\tens(Y\pti\tens X\pti) \xra\sim (X\tens(Y\tens Y\pti))\tens
X\pti \to \\
\xra{X\tens\ev\tens X\pti} (X\tens I)\tens X\pti \xra{r_x\tens X\pti}
X\tens X\pti \xra\ev I,
\end{multline*}
\label{S12.8}\label{p12.8}%
which induces an isomorphism $j_{+X,Y}: Y\pti\tens X\pti \to (X\tens
Y)\pti$, such that the above pairing coincides with
\[ (X\tens Y)\tens(Y\pti\tens X\pti) \xra{1\tens j_+} (X\tens Y)\tens
(X\tens Y)\pti \xra\ev I .\]
The equation
\begin{multline*}
\coev_{X\tens Y} = \bigl( I \xra{\coev_Y} Y\pti\tens Y \simeq
Y\pti\tens I\tens Y \xra{1\tens\coev_X\tens1} \\
Y\pti\tens X\pti\tens X\tens Y \xra{j_+\tens1}
(X\tens Y)\pti\tens(X\tens Y) \bigr)
\end{multline*}
also holds. Similarly, there is an isomorphism $j_{-X,Y}: \lpti Y\tens\lpti
X \to \lpti(X\tens Y)$.

There are canonical isomorphisms
\[ X\to \lpti(X\pti) ,\qquad X \to (\lpti X)\pti .\]
To simplify notations we assume the functors $\text-\pti$ and
$\lpti\text-$ inverse to each other (we can always achieve this replacing
the category by an equivalent one). We shall denote the iterated duals by
$X^{(n\vee)} = X^{\vee\dots\vee}$ ($n$ times) and
$X^{(-n\vee)} = {}^{\vee\dots\vee} X$ ($n$ times) for $n\ge0$.

The graphical notation for the braiding and its inverse is
\[ c = ( c_{X,Y} : X\tens Y \to Y\tens X ) =
\raisebox{-6mm}[8mm][6mm]{
\unitlength 0.70mm
\linethickness{0.4pt}
\begin{picture}(20,25)
\put(0,0){\makebox(0,0)[cb]{$Y$}}
\put(20,0){\makebox(0,0)[cb]{$X$}}
\put(0,20){\makebox(0,0)[ct]{$X$}}
\put(20,20){\makebox(0,0)[ct]{$Y$}}
\put(20,15){\line(-2,-1){20}}
\put(20,5){\line(-2,1){8}}
\put(-0,15){\line(2,-1){8}}
\end{picture}
} \quad, \qquad
c^{-1}  =
\raisebox{-6mm}[8mm][6mm]{
\unitlength 0.70mm
\linethickness{0.4pt}
\begin{picture}(20,25)
\put(0,0){\makebox(0,0)[cb]{$X$}}
\put(20,0){\makebox(0,0)[cb]{$Y$}}
\put(0,20){\makebox(0,0)[ct]{$Y$}}
\put(20,20){\makebox(0,0)[ct]{$X$}}
\put(20,5){\line(-2,1){20}}
\put(-0,5){\line(2,1){8}}
\put(20,15){\line(-2,-1){8}}
\end{picture}
} \quad. \]

In a rigid monoidal braided category there are functorial isomorphisms
\[
\unitlength 0.70mm
\linethickness{0.4pt}
\begin{picture}(146.33,37)
\put(23,18){\oval(10,10)[r]}
\put(23,20){\oval(6,6)[lt]}
\put(11,35){\makebox(0,0)[cc]{$X$}}
\put(11,1){\makebox(0,0)[cc]{$X\lpti\lpti$}}
\put(1,18){\makebox(0,0)[cc]{$u_1^2\ =$}}
\put(61,18){\oval(10,10)[r]}
\put(61,20){\oval(6,6)[lt]}
\put(49,35){\makebox(0,0)[cc]{$X$}}
\put(49,1){\makebox(0,0)[cc]{$X\lpti\lpti$}}
\put(39,18){\makebox(0,0)[cc]{, $u_{-1}^2\ =$}}
\put(92,20){\oval(6,6)[rt]}
\put(106,35){\makebox(0,0)[cc]{$X$}}
\put(106,1){\makebox(0,0)[cc]{$\lpti\lpti X$}}
\put(76,18){\makebox(0,0)[cc]{, $u_1^{-2}=$}}
\put(132,20){\oval(6,6)[rt]}
\put(146,35){\makebox(0,0)[cc]{$X$}}
\put(146,1){\makebox(0,0)[cc]{$\lpti\lpti X$}}
\put(116,18){\makebox(0,0)[cc]{, $u_{-1}^{-2}=$}}
\put(20,20){\line(-2,-3){9.33}}
\put(11,30){\line(2,-3){6.67}}
\put(58,16){\line(-2,3){9.33}}
\put(49,6){\line(3,5){6}}
\put(106,6){\line(-4,5){8}}
\put(135,20){\line(4,-5){11.33}}
\put(146,30){\line(-4,-5){8}}
\put(23,16.50){\oval(6,7)[lb]}
\put(61.50,16){\oval(7,6)[lb]}
\put(92,18){\oval(12,10)[l]}
\put(132,18){\oval(12,10)[l]}
\put(91.50,16){\oval(7,6)[rb]}
\put(95,16){\line(4,5){11.20}}
\put(131.50,16){\oval(7,6)[rb]}
\end{picture}
\]

\subsection{Monoidal abelian categories}\label{S13}
Let $\cv = (\cv, \tens, a, I, r, l)$ be a monoidal abelian category.
Then it has a canonical \kd-linear structure, where $\kk = \End_\cv I$
is a commutative ring, $I$ is the unit object. For $\lambda\in \End_\cv
I$, $f\in\Hom_\cv(X,Y)$ the morphism $\lambda f$ is defined as $(X
\xra{l^{-1}} I\tens X \xra{\lambda\tens f} I\tens Y \xra{l} Y)$.  We
assume that $\cv$ is a rigid monoidal abelian category with length and
the unit object $I$ such that $\kk = \End I$ is a perfect field.  Then
the object $I$ is simple.

The tensor product functor $X\tens\text-$ (resp. $\text-\tens X$) has a
right adjoint $X\pti\tens\text-$ (resp. $\text-\tens\lpti X$) and a
left adjoint $\lpti X\tens\text-$ (resp. $\text-\tens X\pti$). Therefore,
the functor $\tens: \cv\times\cv \to \cv$ is exact in each variable and it
is \kd-bilinear by the choice of \kd-linear structure. By
\propref{P11A} there exists a \kd-linear exact functor $\ot:
\cv\tens\cv\to \cv$ called the diagonal restriction functor, such that
$\tens$ is isomorphic to the composite
\[ \cv\times\cv \xra\odot \cv\tens\cv \xra\ot \cv .\]
The functors $X\mapsto X\pti$, $X\mapsto \lpti X$, quasiinverse to each
other, are also exact.

\begin{proposition}\label{P14}
The functor $\ot: \cv\tens\cv \to \cv$ is faithful.
\end{proposition}

\subsection{System of notations}
The functors $\odot$ and $\ot$ have their analogues acting between
categories $\cv^{\tens n}$ and such. There are several isomorphic functors
of that kind.  We denote their common value by subindices. So, applying the
functor
\[ \cv^{\tens n_1}\times \cv^{\tens n_2}\times \cv^{\tens n_3}\times\dots
\to \cv^{\tens p} \]
to $(A,B,C,\dots)$ we get $A_{i_1\dots i_{n_1}} \tens B_{j_1\dots j_{n_2}}
\tens C_{k_1\dots k_{n_3}} \tens\dots$, where the sign $\tens$ may be
changed to $\odot$ or $\ot$ for purely aesthetical reasons and contains no
additional information. Everything is encoded in terms of indices, which
are all distinct and are taken from the set
\begin{multline*}
\{i_1,\dots,i_{n_1},j_1,\dots,j_{n_2},k_1,\dots,k_{n_3}, \dots\} \\
= \{1',1'',1''',\dots,1^{m_1},2',2'',2''',\dots,2^{m_2},
3',\dots,3^{m_3},\dots,p',\dots,p^{m_p}\} .
\end{multline*}
The number means the tensorand in $\cv^{\tens p}$ to which the present
tensorand goes, and the superscripts determine the order in which several
terms are tensored to get an object of $\cv$ -- tensorand from $\cv^{\tens
p}$.  Another way to describe such functor is to give a permutation from
$\SSS_n$ and two partitions $n = n_1+n_2+n_3+ \dots = m_1+m_2+ \dots +m_p$.

\begin{examples}
If $X,Y\in\cv$, $C\in\vtv$ then $X_{1'}\tens Y_{1''}$ denotes the usual
tensor product of $X$ and $Y$, $X_{1''}\tens Y_{1'}$ is $Y\tens X$.
Similarly, $C_{1'1''}$ is $\ot C$ and $C_{1''1'}$ is $\ot PC$, where $P:
\vtv \to \vtv$ is the permutation functor, $P(X\odot Y) = Y\odot X$. The
three objects $X_1\odot C_{23}$, $X_2\odot C_{13}$, $C_{12}\odot X_3$ of
$\cv^{\tens3}$ differ by the place where $X$ goes. Applying
$\ot\tens\Id_\cv :  \cv^{\tens3} \to \cv^{\tens2}$ to $X_1\odot C_{23}$ we
get $X_{1'}\tens C_{1''2}$. Applying $\Id_\cv\tens\ot : \cv^{\tens3} \to
\cv^{\tens2}$ to $C_{12}\odot X_3$ we get $C_{12'}\tens X_{2''}$. The
functor $\Id_\cv\tens \ot\tens\Id_\cv: \cv^{\tens4} \to \cv^{\tens3}$
(tensoring the second and the third places), applied to $C\odot C$, gives
$C_{12'}\tens C_{2''3}$.
\end{examples}

To use the same system of notations for morphisms we should write indices
for both source and target. However, we shall abbreviate the notation using
only one set of indices, either for the source, or for the target. Thus,
instead of
\begin{gather*}
{\ }_{1'1''}\e_1 : C_{1'1''} \to I_1 ,\\
\ _{12'}C_{12'} \tens {}_{2''2'''}\e_{2''} : C_{12'}\tens C_{2''2'''} \to
C_{12'}\tens I_{2''}
\end{gather*}
we write $\e: C_{1'1''}\to I$ and $C_{12'}\tens\e_{2''}$. Instead of
\begin{gather*}
{}_{123}\Delta_{12'2''3} :  C_{13}\odot I_2 \to C_{12'}\tens C_{2''3} ,\\
{}_{124}\Delta_{12'2''4}\odot I_3 : C_{14}\odot I_2\odot I_3 \to
C_{12'}\tens C_{2''4}\odot I_3
\end{gather*}
we write $\Delta_{123}$ and $\Delta_{124}\odot I_3$.

Sometimes we simplify $X_{1'}$, $X_{1''}$, $X_{1'''}$, $X_{1^4}, \dots$
to $X^1$, $X^2$, $X^3$, $X^4, \dots$.

\subsection{Monoidal structures on $\vtv$}\label{S15.4}
Let $\cv$ be a \kd-linear abelian monoidal category with
length.  Then $\vtv$ has a monoidal structure as well. In fact, there are
four monoidal structures, since we may choose between $(\cv,\tens)$ and
$(\cv,\tens_\op)$ in both tensorands. The main monoidal structure, which we
fix from now on, is:
\begin{align*}
\barten: (\vtv)\times(\vtv) &\longrightarrow \vtv\\
(A_{12},B_{12}) &\longmapsto A_{1'2''}\tens B_{1''2'}\\
(X\odot Y,V\odot W) &\longmapsto (X\tens V)\odot (W\tens Y) .
\end{align*}

\begin{theorem}\label{T15.4}
If $(\cv,\tens)$ is rigid, then $(\vtv,\barten)$ is rigid as well.
\end{theorem}

\subsection{$\Ind$-objects}
Following Grothendieck and Verdier \cite{GroVer} we consider the category
of $\ind$-objects of a given \kd-linear abelian category $\ca$. Recall that
an $\ind$-object of $\ca$ is a functor $X:  I\to\ca$ from a filtered
partially ordered set $I$, in particular, an arrow $x_{ij} : X_i \to X_j
\in\ca$ is given for $i<j$, $i,j\in I$. The set of morphisms from $X: I
\to\ca$ to $Y: J \to\ca$ is
\[ \limit{I} \colimit{J} \Hom(X_i,Y_j) . \]
We denote the category of $\ind$-objects of $\ca$ by $\ahat$ as a synonym
of standard notation $\Ind(\ca)$. The category $\ahat$ is a \kd-linear
abelian category \cite[Exercise 8.9.9]{GroVer}. Small projective and
inductive limits in $\ahat$ are representable in $\ahat$ (Grothendieck and
Verdier \cite[Propositions 8.9.1 and 8.9.5]{GroVer}).

\begin{theorem}[Grothendieck and Verdier \cite{GroVer} Theorem 8.3.3]
Let $\ca$ be essentially small, then the functor
\begin{align*}
\ahat &\longrightarrow \und\Hom_{\kk,\text{l.e.}} (\ca^\op, \kk\Vect) ,\\
(X:I \to\ca) &\longmapsto (Y \mapsto \colimit{I} \Hom(Y, X_i))
\end{align*}
with values in the category of \kd-linear left exact functors is an
equivalence of categories.
\end{theorem}

If $\ca,\cb$ are \kd-linear abelian essentially small categories, then any
functor $F: \ca \to \cb$ extends to a functor $\hat F: \ahat \to \bhat$, $X
\mapsto F\circ X$. If $F$ is \kd-linear (resp. right exact, resp.  exact),
so is $\hat F$ (Grothendieck and Verdier \cite[Corollary 8.9.8]{GroVer}).
The functor $\hat F: \ahat\to\bhat$ commutes with filtered inductive limits
\cite[Proposition 8.6.3]{GroVer}. If $F$ is right exact, $\hat F$ commutes
with arbitrary inductive limits.

\begin{proposition}
Let any object of $\ca$ have finite length. Then the category $\ahat$ is
equivalent to its full subcategory consisting of functors $X: I \to\ca$
such that $x_{ij} : X_i \to X_j$ is a monomorphism for any pair $i<j$.
\end{proposition}

\begin{remark}
Let $X:I \to\ca$ be in $\ahat$ and let $J\subset I$ be a cofinal set. Then
the $\ind$-object $X' = X\vert_J :J \to \ca$ is isomorphic to $X$ in
$\ahat$.
\end{remark}

Let $X: I \to\ca$ be an $\ind$-object such that $x_{ij} : X_i \to X_j$ are
monomorphisms, $i<j$. We say that $X$ is {\em closed under intersections}
if
\begin{enumerate}
\renewcommand{\labelenumi}{(\alph{enumi})}
\item for any subset $J\subset I$ there is an element $i = \cap J \in I$
such that $i\le j$ for any $j\in J$ and $i$ is the biggest element with
this property;
\item for any subset $J\subset I$ there is a finite subset $J'\subset
J$ such that for any finite $K$, $J'\subset K\subset J$, and any $r\ge K$
the subobject $X_i$ is the intersection of subobjects $X_k$, $k\in K$, in
$X_r$, that is, the canonical morphism
\[ X_i \to \limit{k\in K} (X_k \xra{x_{kr}} X_r) \]
is an isomorphism.
\end{enumerate}

\begin{proposition}
\textup{(a)} Any $\ind$-object is isomorphic in $\ahat$
to an $\ind$-object closed under intersections.

\textup{(b)}
Let $X: I\to\ca$, $Y:J\to\ca$ be $\ind$-objects. Assume that $Y$ is
closed under intersections. Then any morphism $f:X\to Y\in\ahat$ can be
represented by a monotonous map $m:I\to J$ and a family of morphisms
$f_i: X_i \to Y_{m(i)}$.
\end{proposition}

\section{Squared coalgebras}
\subsection{Definitions}
Let $\cv$ be a \kd-linear abelian rigid monoidal category with
 length.

\begin{definition}
A {\em squared coalgebra} $C = (C,\Delta,\e)$ in $\vhat$ is an object
$C\in\vtvhat$ equipped with a comultiplication
$\Delta_{123}: C_{13}\odot I_2 \to C_{12'}\tens C_{2''3}
\in\wh{\cv^{\tens3}}$ and a counit $\e: C_{1'1''} \to I \in\vhat$,
such that coassociativity holds:
\begin{equation}\label{d23a}
\begin{CD}
C_{14}\odot I_2\odot I_3
@>\Delta_{124}\odot I_3>\hphantom{\Delta_{123'}\tens C_{3''4}}>
C_{12'}\tens C_{2''4}\odot I_3 \\
@V\Delta_{134}\odot I_2VV       @VVC_{12'}\tens\Delta_{2''34}V  \\
C_{13'}\tens C_{3''4}\odot I_2 @>\Delta_{123'}\tens C_{3''4}>>
C_{12'}\tens C_{2''3'}\tens C_{3''4}
\end{CD}
\end{equation}
and $\e$ is the counit:
\begin{subequations}
\begin{gather}\label{e23b}
\bigl( C_{12} \xra{\sim} C_{1'2}\tens I_{1''} \xra{\Delta_{1'1''2}}
C_{1'1''}\tens C_{1'''2} \xra{\e\tens C} I_{1'}\tens C_{1''2}
\xra\sim C_{12} \bigr) = \id_C \\
\bigl( C_{12} \xra{\sim} I_{2'}\tens C_{12''} \xra{\Delta_{12'2''}}
C_{12'}\tens C_{2''2'''} \xra{C\tens\e} C_{12'}\tens I_{2''}
\xra\sim C_{12} \bigr) = \id_C \label{e23c}
\end{gather}
\end{subequations}
\end{definition}

Graphical notations partially explain the choice of indices and help to
memorise them. Comultiplication is denoted
\label{p23}
\[
\unitlength 1mm
\linethickness{0.4pt}
\begin{picture}(48.33,19.63)
\put(15,1.08){\makebox(0,0)[cc]{$C$}}
\put(35,1.29){\makebox(0,0)[cc]{$C$}}
\put(9.83,1.58){\makebox(0,0)[cb]{1}}
\put(20,1.58){\makebox(0,0)[cb]{$2'$}}
\put(30,1.58){\makebox(0,0)[cb]{$2''$}}
\put(40,1.58){\makebox(0,0)[cb]{3}}
\put(25,4.95){\oval(10,9.89)[t]}
\put(40,4.95){\line(0,1){10.11}}
\put(10,4.95){\line(0,1){10.11}}
\put(9.83,16.34){\makebox(0,0)[cb]{1}}
\put(40,16.34){\makebox(0,0)[cb]{3}}
\put(25.17,19.63){\makebox(0,0)[ct]{$C$}}
\put(2.50,9.89){\makebox(0,0)[rb]{$\Delta$ $=$}}
\put(48.33,9.89){\makebox(0,0)[cc]{,}}
\end{picture}
\]
the counit is denoted
\[
\unitlength 1mm
\linethickness{0.4pt}
\begin{picture}(25,12.04)
\put(14.94,6.02){\oval(9.89,12.04)[b]}
\put(9.89,7.03){\makebox(0,0)[cb]{$1'$}}
\put(20,7.03){\makebox(0,0)[cb]{$1''$}}
\put(15,10.04){\makebox(0,0)[ct]{$C$}}
\put(5,5.02){\makebox(0,0)[rb]{$\e$ $=$}}
\put(25,5.02){\makebox(0,0)[cc]{.}}
\end{picture}
\]
The coassociativity equation is
\[
\unitlength 0.80mm
\linethickness{0.4pt}
\begin{picture}(140,24.09)
\put(15,0.22){\makebox(0,0)[cb]{$C$}}
\put(35,0.22){\makebox(0,0)[cb]{$C$}}
\put(55,0.22){\makebox(0,0)[cb]{$C$}}
\put(10.16,0.65){\makebox(0,0)[cb]{1}}
\put(20.33,0.65){\makebox(0,0)[cb]{$2'$}}
\put(30.33,0.65){\makebox(0,0)[cb]{$2''$}}
\put(40.33,0.65){\makebox(0,0)[cb]{$3'$}}
\put(50.33,0.65){\makebox(0,0)[cb]{$3''$}}
\put(60.66,0.65){\makebox(0,0)[cb]{4}}
\put(44.92,4.95){\oval(10.17,9.89)[t]}
\put(25,5.05){\oval(10,20)[t]}
\put(10,4.95){\line(0,1){15.05}}
\put(60.17,4.73){\line(0,1){15.27}}
\put(10,21.51){\makebox(0,0)[cb]{1}}
\put(60.17,21.51){\makebox(0,0)[cb]{4}}
\put(35.17,24.09){\makebox(0,0)[ct]{$C$}}
\put(85,0.22){\makebox(0,0)[cb]{$C$}}
\put(105,0.22){\makebox(0,0)[cb]{$C$}}
\put(125,0.22){\makebox(0,0)[cb]{$C$}}
\put(80.16,0.65){\makebox(0,0)[cb]{1}}
\put(90.33,0.65){\makebox(0,0)[cb]{$2'$}}
\put(100.33,0.65){\makebox(0,0)[cb]{$2''$}}
\put(110.33,0.65){\makebox(0,0)[cb]{$3'$}}
\put(120.33,0.65){\makebox(0,0)[cb]{$3''$}}
\put(130.66,0.65){\makebox(0,0)[cb]{4}}
\put(95.01,4.95){\oval(10.17,9.89)[t]}
\put(114.91,5.05){\oval(10,20)[t]}
\put(80,4.95){\line(0,1){15.05}}
\put(130.17,4.73){\line(0,1){15.27}}
\put(80,21.51){\makebox(0,0)[cb]{1}}
\put(130.17,21.51){\makebox(0,0)[cb]{4}}
\put(105.17,24.09){\makebox(0,0)[ct]{$C$}}
\put(70.33,12.04){\makebox(0,0)[cc]{$=$}}
\put(140,12.04){\makebox(0,0)[cc]{,}}
\end{picture}
\]
the equations for the counit are
\[
\raisebox{-11mm}{
\unitlength 0.80mm
\linethickness{0.4pt}
\begin{picture}(123.50,29.03)
\put(14.92,14.95){\oval(10.17,10.11)[b]}
\put(25,14.95){\oval(10,10.11)[t]}
\put(30,15.05){\line(0,-1){10.11}}
\put(9.83,15.05){\line(0,1){9.89}}
\put(40,24.95){\line(0,-1){20}}
\put(50,15.05){\makebox(0,0)[cc]{$=$}}
\put(6.50,13.98){\makebox(0,0)[cb]{$1'$}}
\put(17,13.98){\makebox(0,0)[cb]{$1''$}}
\put(33,13.98){\makebox(0,0)[cb]{$1'''$}}
\put(30,4.09){\makebox(0,0)[ct]{1}}
\put(40,4.09){\makebox(0,0)[ct]{2}}
\put(35,0.22){\makebox(0,0)[cb]{$C$}}
\put(9.83,26.02){\makebox(0,0)[cb]{1}}
\put(40,26.02){\makebox(0,0)[cb]{2}}
\put(25.17,29.03){\makebox(0,0)[ct]{$C$}}
\put(60,4.95){\line(0,1){20}}
\put(70,24.95){\line(0,-1){20}}
\put(60,4.09){\makebox(0,0)[ct]{1}}
\put(70,4.09){\makebox(0,0)[ct]{2}}
\put(65,0.22){\makebox(0,0)[cb]{$C$}}
\put(60,26.02){\makebox(0,0)[cb]{1}}
\put(70,26.02){\makebox(0,0)[cb]{2}}
\put(65,29.03){\makebox(0,0)[ct]{$C$}}
\put(115.08,14.95){\oval(10.17,10.11)[b]}
\put(105,14.95){\oval(10,10.11)[t]}
\put(100,15.05){\line(0,-1){10.11}}
\put(120.17,15.05){\line(0,1){9.89}}
\put(90,24.95){\line(0,-1){20}}
\put(80,15.05){\makebox(0,0)[cc]{$=$}}
\put(123.50,13.98){\makebox(0,0)[cb]{$2'''$}}
\put(113,13.98){\makebox(0,0)[cb]{$2''$}}
\put(97,13.98){\makebox(0,0)[cb]{$2'$}}
\put(100,4.09){\makebox(0,0)[ct]{2}}
\put(90,4.09){\makebox(0,0)[ct]{1}}
\put(95,0.22){\makebox(0,0)[cb]{$C$}}
\put(120.17,26.02){\makebox(0,0)[cb]{2}}
\put(90,26.02){\makebox(0,0)[cb]{1}}
\put(104.83,29.03){\makebox(0,0)[ct]{$C$}}
\end{picture}
} \quad .
\]

\begin{definition}
A {\em squared coalgebra homomorphism}
$(C,\Delta,\e_C) \to (D,\Delta,\e_D)$ is a morphism
$f:C\to D \in\Mor\widehat{\cv\tens\cv}$ such that
\begin{equation*}
\begin{CD}
C_{13}\odot I_2 @>\Delta_{123}>> C_{12'}\tens C_{2''3} \\
@Vf_{13}\odot I_2VV             @VVf_{12'}\tens f_{2''3}V \\
D_{13}\odot I_2 @>\Delta_{123}>> D_{12'}\tens D_{2''3}
\end{CD}
\end{equation*}
commutes and equation
\begin{equation*}
(\ot C \xra{\ot f} \ot D \xra{\e_D} I) = \e_C
\end{equation*}
holds.
\end{definition}

Squared coalgebras form a category denoted $\Coalgsq(\vhat)$. Its full
subcategory consisting of squared coalgebras, which are objects of
$\vtv$, is denoted $\Coalgsq(\cv)$.

\subsection{Comodules}
\begin{definition}
A {\em left comodule} $X\in\vhat$ over a squared coalgebra $C$ is an
object $X$ of $\vhat$ equipped with the coaction
\[ \delta = \delta_X : X_1\odot I_2 \to C_{12'}\tens X_{2''}
\in \widehat{\cv\tens\cv} , \]
which is coassociative:
\begin{equation}\label{d28a}
\begin{CD}
X_1\odot I_2\odot I_3
@>\delta_{12}\odot I_3>\hphantom{\Delta_{123'}\tens X_{3''}}>
C_{12'}\tens X_{2''}\odot I_3 \\
@V\delta_{13}\odot I_2VV        @VVC_{12'}\tens\delta_{2''3}V \\
C_{13'}\odot I_2\tens X_{3''} @>\Delta_{123'}\tens X_{3''}>>
C_{12'}\tens C_{2''3'}\tens X_{3''}
\end{CD}
\end{equation}
and counital:
\begin{equation}\label{e28b}
\left(X \xra\sim X\tens I \xra{\ot\delta} (\ot C)\tens X \xra{\e\tens X}
I\tens X \xra\sim X\right) = \id_X .
\end{equation}
\end{definition}

Graphical notation for the coaction is
\label{p29}
\[
\unitlength 1mm
\linethickness{0.4pt}
\begin{picture}(35,18.06)
\put(15,1.08){\makebox(0,0)[cc]{$C$}}
\put(35,1.29){\makebox(0,0)[cc]{$X$}}
\put(9.74,1.57){\makebox(0,0)[cb]{1}}
\put(19.91,1.57){\makebox(0,0)[cb]{$2'$}}
\put(29.91,1.57){\makebox(0,0)[cb]{$2''$}}
\put(25,4.95){\oval(10,9.89)[t]}
\put(10,4.95){\line(0,1){10.11}}
\put(9.83,16.34){\makebox(0,0)[cb]{1}}
\put(14.84,18.06){\makebox(0,0)[ct]{$X$}}
\put(2.50,9.89){\makebox(0,0)[rb]{$\delta$ $=$}}
\end{picture}
\]
and it will be explained later. Coassociativity takes the form
\[
\unitlength 0.80mm
\linethickness{0.4pt}
\begin{picture}(115.33,24.09)
\put(15,0.22){\makebox(0,0)[cb]{$C$}}
\put(35,0.22){\makebox(0,0)[cb]{$C$}}
\put(55,0.22){\makebox(0,0)[cb]{$X$}}
\put(9.83,0.66){\makebox(0,0)[cb]{1}}
\put(20,0.66){\makebox(0,0)[cb]{$2'$}}
\put(30,0.66){\makebox(0,0)[cb]{$2''$}}
\put(40,0.66){\makebox(0,0)[cb]{$3'$}}
\put(50,0.66){\makebox(0,0)[cb]{$3''$}}
\put(44.92,4.95){\oval(10.17,9.89)[t]}
\put(25,5.05){\oval(10,20)[t]}
\put(10,4.95){\line(0,1){15.05}}
\put(10,21.51){\makebox(0,0)[cb]{1}}
\put(75.33,0.22){\makebox(0,0)[cb]{$C$}}
\put(95.33,0.22){\makebox(0,0)[cb]{$C$}}
\put(115.33,0.22){\makebox(0,0)[cb]{$X$}}
\put(70.16,0.66){\makebox(0,0)[cb]{1}}
\put(80.33,0.66){\makebox(0,0)[cb]{$2'$}}
\put(90.33,0.66){\makebox(0,0)[cb]{$2''$}}
\put(100.33,0.66){\makebox(0,0)[cb]{$3'$}}
\put(110.33,0.66){\makebox(0,0)[cb]{$3''$}}
\put(85.34,4.95){\oval(10.17,9.89)[t]}
\put(105.24,5.05){\oval(10,20)[t]}
\put(70.33,4.95){\line(0,1){15.05}}
\put(70.33,21.51){\makebox(0,0)[cb]{1}}
\put(60.66,12.04){\makebox(0,0)[cc]{$=$}}
\put(74.67,24.09){\makebox(0,0)[ct]{$X$}}
\put(14.67,23.09){\makebox(0,0)[ct]{$X$}}
\end{picture}
\]
and counitality is
\[
\unitlength 0.80mm
\linethickness{0.4pt}
\begin{picture}(69.33,29.03)
\put(14.92,14.95){\oval(10.17,10.11)[b]}
\put(25,14.95){\oval(10,10.11)[t]}
\put(30,15.05){\line(0,-1){10.11}}
\put(9.83,15.05){\line(0,1){9.89}}
\put(43.33,15.05){\makebox(0,0)[cc]{$=$}}
\put(6.50,13.98){\makebox(0,0)[cb]{$1'$}}
\put(17,13.98){\makebox(0,0)[cb]{$1''$}}
\put(33,13.98){\makebox(0,0)[cb]{$1'''$}}
\put(30,4.09){\makebox(0,0)[ct]{1}}
\put(35,0.22){\makebox(0,0)[cb]{$X$}}
\put(9.83,26.02){\makebox(0,0)[cb]{1}}
\put(14.84,28.60){\makebox(0,0)[ct]{$X$}}
\put(53.33,4.95){\line(0,1){20}}
\put(53.33,4.09){\makebox(0,0)[ct]{1}}
\put(58.33,0.22){\makebox(0,0)[cb]{$X$}}
\put(53.33,26.02){\makebox(0,0)[cb]{1}}
\put(58.33,29.03){\makebox(0,0)[ct]{$X$}}
\put(69.33,14.62){\makebox(0,0)[cc]{.}}
\end{picture}
\]
The definition above should be generalised. Let $\ca$ be a \kd-linear
abelian category with length.

\begin{definition}
A {\em left comodule} $X\in\vtahat$ over a squared coalgebra
$C\in\Coalgsq(\vhat)$ is an object $X = X_{10} \in\vtahat$ equipped
with the coaction
\[ \delta = \delta_X : X_{10}\odot I_2 \to C_{12'}\tens X_{2''0}
\in (\cv\tens\vta)\sphat \]
such that coassociativity \eqref{d28a} in $(\vto3\tens\ca)\sphat$ and
counitality \eqref{e28b} in $\vtahat$ hold.
\end{definition}

When $\ca = \kk\vect$ this reduces to the previous definition.

\begin{definition}
A morphism of left \Cd-comodules $(X,\delta_X)\to(Y,\delta_Y)$ in
$\vhat$ (resp. in $\vtahat$) is $f:X\to Y \in\vhat$ (resp.
$f: X_{10} \to Y_{10} \in\vtahat$) such that
\begin{equation*}
\begin{CD}
X_1\odot I_2 @>\delta_X>> C_{12'}\tens X_{2''} \\
@Vf_1\odot I_2VV        @VVC_{12'}\tens f_{2''}V \\
Y_1\odot I_2 @>\delta_Y>> C_{12'}\tens Y_{2''}
\end{CD}
\end{equation*}
(resp. the same diagram in $(\vtv\tens\ca)\sphat$\ ) commutes.
\end{definition}

Left \Cd-comodules form a category, which is denoted ${}^C\vhat$ (resp.
${}^C\vtahat$). It has a full subcategory ${}^C\cv$ (resp.
${}^C(\vta)$) formed by objects from $\cv$ (resp. $\vta$).

It is easy to show that if $f:X\to Y$ is a morphism of comodules in
${}^C\vhat$ (resp. ${}^C\vtahat$, ${}^C\cv$, ${}^C(\vta)$), and
\[ \Ker f \xra{\kernel f} X \xra{\coim f} \Coim f \simeq \Im f
\xra{\im f} Y \xra{\coker f} \Coker f \]
is its canonical decomposition in $\vhat$ (resp. $\vtahat$,  $\cv$,
$\vta$), then the objects $\Ker f$, $\Coim f$, $\Im f$, $\Coker f$
have unique \Cd-comodule structure such that the morphisms above are
morphisms of comodules. It follows:

\begin{proposition}\label{P30}
The category ${}^C\vhat$ (resp. ${}^C\vtahat$, ${}^C\cv$,
${}^C(\vta)$) is \kd-linear and abelian and the underlying functor
$\cu: {}^C\vhat \to \vhat$ (resp. $\cu: {}^C\vtahat \to \vtahat$, $\cu:
{}^C\cv \to \cv$, $\cu: {}^C(\vta) \to \vta$) is exact and
faithful.
\end{proposition}

Any \kd-linear exact functor $F:\ca\to\cb$ induces a \kd-linear exact
functor $\cv\tens F: {}^C\vtahat \to {}^C\wh{\cv\tens\cb}$.

\begin{example}
$(C,\Delta)$ is a \Cd-comodule from ${}^C\vtvhat$ (resp. $\vtvhat^C$). It
is called the left 
regular comodule.
\end{example}

\subsection{The fundamental theorem on coalgebras}
\begin{theorem}\label{T31}
Any comodule from ${}^C\vhat$ (resp. ${}^C\vtahat$) is a union, i.e.
filtered inductive limit, of its subcomodules from ${}^C\cv$ (resp.
${}^C(\vta)$).
\end{theorem}

\begin{corollary}
$\wh{{}^C\cv} \simeq {}^C\vhat$ and
$\bigl({}^C(\vta)\bigr)\sphat \simeq {}^C\vtahat$.
\end{corollary}

\begin{theorem}[fundamental theorem on coalgebras]\label{T34}
A squared coalgebra $C\in \Coalgsq(\vhat)$ is a filtered inductive limit of
its subcoalgebras from $\Coalgsq(\cv)$.
\end{theorem}

Let $\ca$ be a \kd-linear abelian category with length
and let $C$ be a squared coalgebra in $\cv$.

\begin{theorem}\label{T34.1}
The functor $\Psi: ({}^C\cv)\tens\ca \to {}^C(\vta)$ induced by
$(^C\cv)\times\ca \to {}^C(\vta)$, $(X,M) \mapsto X\odot M$,
is an equivalence.
\end{theorem}

\subsection{The category of fibre functors}\label{S35}
Let $\cv$ be a \kd-linear abelian category. Extending the definition of
Saavedra~\cite{SaaRiv} let us call a {\em fibre functor} to $\cv$ a
\kd-linear exact faithful functor $a:\ca\to\cv$, where $\ca$ is a
\kd-linear abelian essentially small category. Following
Schauenburg~\cite{Sch:hopf1} we define the category of fibre functors.

\begin{definition}\label{D35}
Let the category $\fa = \fa(\cv)$ have fibre functors to $\cv$ as objects
and let morphisms from $a:\ca\to\cv$ to $b:\cb\to\cv$ be equivalence
classes of pairs $(F,\phi)$, where $F:\ca\to\cb$ is a functor and $\phi: a
\xra\sim bF$ is a functorial isomorphism. Two such pairs $(F,\phi)$ and
$(G,\gamma)$ are equivalent if there is a functorial isomorphism
$\zeta:F\to G$ such that
\begin{equation}\label{e35}
\gamma = \bigl( a \xra\phi bF \xra{b\zeta} bG \bigr) .
\end{equation}
The composite of two morphisms represented by $(F,\phi)$ and $(G,\gamma)$
is represented by $(GF,\gamma F\circ\phi)$. Clearly, $\Hom_\fa(a,b)$ is a
set.
\end{definition}

Now let $\cv$ be a \kd-linear abelian rigid monoidal category with
 length.

There is a functor $\Phi: \Coalgsq(\vhat) \to \fa(\cv)$, $C \mapsto (\cu:
{}^C\cv \to \cv)$, where $\cu$ is the underlying functor. To a squared
coalgebra morphism $f:C\to D$ corresponds the equivalence class of the pair
$(F,\phi)$, where
\[ F(X,\delta_X) = (X, X_1\odot I_2 \xra{\delta_X} C_{12'}\tens X_{2''}
\xra{f\tens X} D_{12'}\tens X_{2''}) \]
and $\phi$ is the identity automorphism $\cu \to \cu$.

Our primary goal is to show that the functor $\Phi$ is an equivalence.

\subsubsection{The coend}\label{S37}
Let $\cc$ be a \kd-linear abelian category with length. The coend $C$
of a bifunctor $B: \cp\times\cp^\op \to \cc$ is defined in \cite{Mac:work}
as an object of $\chat$ which is the inductive limit of the diagram
\[ B(X,X) \xla{B(X,f)} B(X,Y) \xra{B(f,Y)} B(Y,Y) ,\]
where $f:X\to Y$ runs over $\Mor\cp$. That is, $C$ is equipped with a
morphism $i_X: B(X,X) \to C \in \chat$ for each $X\in\Ob\cp$, the diagram
\[ \begin{CD}
B(X,Y) @>B(f,Y)>> B(Y,Y) \\
@VB(X,f)VV  @VVi_YV \\
B(X,X) @>i_X>\hphantom{B(f,Y)}> C
\end{CD} \]
is commutative for any $f:X \to Y\in \Mor\cp$, and $C$ is universal between
such objects. If $\cp$ is small, we can say that the sequence
\begin{equation*}
\bigoplus_{f:X\to Y\in\Mor\cp} B(X,Y) \xra{B(X,f)-B(f,Y)}
\bigoplus_{X\in\Ob\cp} B(X,X) \xra{\oplus i_X} C \to 0
\end{equation*}
is exact. So in this case the coend exists. More generally, it exists for
essentially small $\cp$. It is denoted $C= \int^{X\in\cp} B(X,X)$.

Let us consider the particular case. Let $p:\cp\to\cv$ be a functor from
an essentially small category $\cp$, let $\cc = \vtv$ and let
$B:\cp\times\cp^\op \to\cc$, $B(X,Y) = pX\odot(pY)\pti$.
The coend is denoted
\begin{equation}\label{e38a}
 C = \int^{X\in\cp} pX\odot(pX)\pti .
\end{equation}

The object $M\odot M\pti \in\vtv$ has a canonical squared coalgebra
structure for any $M\in\cv$. Namely,
\begin{gather*}
\Delta_{123} = M\odot\coev\odot M\pti : M\odot I\odot M\pti \to M\odot
M\pti\tens M\odot M\pti ,\\
\e = \ev : M\tens M\pti \to I.
\end{gather*}
Compare these formulae with the graphical notations on page~\pageref{p23}.
The object $M$ has a canonical structure of a left $M\odot M\pti$-comodule,
namely,
\label{S38}
\[ \delta = M\odot\coev : M_1\odot I_2 \to M_1\odot M\pti_{2'} \tens
M_{2''} .\]
Compare with the graphical notations on page~\pageref{p29}. In particular,
this holds for $M = pX$.

\begin{proposition}\label{P38}
If the coend \eqref{e38a} exists, it has a unique squared coalgebra
structure such that the structure morphisms $i_X: pX\odot (pX)\pti \to C$
are coalgebra morphisms.
\end{proposition}

Since $pX$ is a $pX\odot(pX)\pti$-comodule, it is a \Cd-comodule as well
for any $X\in \Ob\cp$.

\begin{proposition}\label{P40}
If the coend $C$ exists, the map $\Ob\cp \to \Ob{}^C\cv$, $X \mapsto
(pX,\delta)$ extends to the functor $F:\cp \to {}^C\cv$, $f\mapsto pf$
such that $p = \bigl( \cp \xra{F} {}^C\cv \xra\cu \cv \bigr)$.
\end{proposition}

\subsection{Reconstruction theorems}
Now we come back to the case of a fibre functor $a:\ca\to\cv$. Assume that
$\cp\subset\ca$ is a full subcategory equivalent to $\ca$ and $\co =
\Ob\cp$ is a set. Denote $p= a\vert_\cp : \cp \to\cv$. Notice, that
coend~\eqref{e38a} serves as the coend of the bifunctor $\ca^\op\times\ca
\to \vtv$, $(X,Y) \mapsto aX\odot(aY)\pti$ as well.

\begin{theorem}\label{T49}
The functor $F:\ca \to {}^C\cv$ is an equivalence of categories.
\end{theorem}

In particular case $\cv = \kk\vect$ this theorem was proved by
Saavedra~\cite{SaaRiv} (see also Schauenburg~\cite{Sch:hopf1}).

\begin{theorem}
\textup{(a)}
The map $(a:\ca\to\cv) \longmapsto C_a$ constructed in \secref{S37}
extends to a functor
\[ \Psi : \fa(\cv) \to \Coalgsq(\vhat) .\]

\textup{(b)}
The functor $F_a: \ca \to \Phi(\Psi(\ca)) = {}^{C_a}\cv$ constructed in
\thmref{T49} together with $\id: a\to \cu\circ F_a$ gives an isomorphism of
functors
\[ F: \Id_{\fa(\cv)} \xra\sim \Phi\Psi .\]

\label{T54}\textup{(c)}
The functors $\Phi: \Coalgsq(\vhat) \to \fa(\cv)$, $C \mapsto {}^C\cv$, and
$\Psi: \fa(\cv) \to \Coalgsq(\vhat)$, $a\to C_a$, are equivalences,
quasiinverse to each other.
\end{theorem}

In the case of $\cv=\kk\vect$ this theorem was proved by
Schauenburg~\cite{Sch:hopf1}.

Constructing another adjunction $\Psi\Phi \xra\sim \Id$, is also of
practical interest.

\begin{proposition}\label{P56A}
\textup{(a)}
Let $C$ be a squared coalgebra in $\vhat$ and let $X\in\cv$. The structures
of \Cd-comodules in $X$ and squared coalgebra homomorphisms $X\odot X\pti
\to C$ are in bijective correspondence.

\textup{(b)}
If $\delta: X_1\odot I_2 \to C_{12'}\tens X_{2''}$ is a comodule structure,
then
\begin{multline}\label{e56a}
\ii = \ii_X = \bigl( X_1\odot X_2\pti \xra\sim X_1\odot I_{2'}\tens
X_{2''}\pti \xra{\delta_{12'}\tens X_{2''}\pti} \\
C_{12'}\tens X_{2''}\tens X_{2'''}\pti \xra{C_{12'}\tens \ev_{2''}}
C_{12'}\tens I_{2''} \xra\simeq C_{12} \bigr)
\end{multline}
is a homomorphism of squared coalgebras.

\textup{(c)} If $\ii_X: X_1\odot X_2\pti \to C_{12}$
is a squared coalgebra homomorphism, then
\[ \delta = \bigl( X_1\odot I_2 \xra{X_1\odot\coev_2} X_1\odot
X_{2'}\pti\tens X_{2''} \xra{\ii_{12'}\tens X_{2''}} C_{12'}\tens X_{2''}
\bigr) \]
is a \Cd-comodule structure on $X$.

\label{P56B}\textup{(d)}
There is a unique squared coalgebra homomorphism $h_C: C'
\overset{\text{def}}= \Psi(\Phi C) \to C$ such that for any $X\in {}^C\cv$
the composite coalgebra morphism $\ii_X: X\odot X\pti \xra{i_X} C'
\xra{h_C} C$ is the canonical \eqref{e56a}.

\label{P57}\textup{(e)}
The family $h_C: \Psi(\Phi C) \to C$ gives an isomorphism of functors $h:
\Psi\Phi \to \Id_{\Coalgsq(\vhat)}$.
\end{proposition}

\begin{corollary}\label{C58}
Any squared coalgebra $C$ is the union of the images of the canonical
morphisms~\eqref{e56a} for all comodules $X$.
\end{corollary}

This is a detailed form of the fundamental theorem on coalgebras \ref{T34}.

\subsection{Comodules over ordinary coalgebras}
Let $C$ be a squared coalgebra in $\cv$. Then $\bar C = \ot C \in\cv$ is an
ordinary coalgebra in $\cv$. How big is the category of \Cbar-comodules in
comparison with \Cd-comodules?

First of all, any \Cbar-comodule is a union of its rigid subcomodules,
${}^{\bar C}\vhat = \wh{{}^{\bar C}\cv}$. This is a general fact. It
follows from the embedding $\delta: X \hookrightarrow \bar C\tens X$ of
\Cbar-comodules and the fundamental theorem on coalgebras.

\begin{theorem}\label{T58.2}
The functor
\begin{gather*}
{}^C\cv\times\cv \xra\tens {}^{\bar C}\cv, \qquad (X,\delta_X)\times Y
\mapsto (X\tens Y, \bar\delta_X\tens Y) \\
\bar\delta_X = \bigl( X \simeq X\tens I \xra{\ot\delta_X} C_{1'1''}\tens
X_{1'''} = \bar C\tens X \bigr)
\end{gather*}
makes ${}^{\bar C}\cv$ into ${}^C\cv\tens\cv$.
\end{theorem}

\begin{corollary}
Let $\ca$ be a \kd-linear abelian essentially small category,
let $\omega: \ca \to \cv$ be a \kd-linear exact faithful functor,
and $\bar C = \int^{X\in\ca} \omega X\tens(\omega X)\pti$. Then the
category of \Cbar-comodules ${}^{\bar C}\cv$ is equivalent to
$\ca\tens\cv$.
\end{corollary}

This is closely related with results of Pareigis~\cite{Par:coend}.

\section{Squared bicoalgebras}
Before we discuss bicoalgebras, let us consider the operation of tensor
product of squared coalgebras.  Notice that a braiding in $\cv$ is not
required.

\subsection{Tensor product of squared coalgebras}
Recall that $\cv\tens\cv$ has a rigid monoidal structure $(A\barten B)_{12}
= A_{1'2''}\tens B_{1''2'}$ (\thmref{T15.4}).

\begin{proposition}\label{P66}
Let $A,B \in \Coalgsq(\vhat)$. The tensor product $A\barten B$ has a
coalgebra structure
\begin{multline}\label{e66a}
\Delta^{A\barten B} : A_{1'3''}\tens B_{1''3'}\odot I_2
\xra{A_{1'3''}\tens\Delta_{1''23'}} A_{1'3''}\tens B_{1''2'}\tens B_{2''3'}
\simeq \\
A_{1'3''}\tens B_{1''2'}\tens I_{2''}\tens B_{2'''3'}
\xra{\Delta_{1'2''3''}\tens B_{1''2'}\tens B_{2'''3'}}
A_{1'2''}\tens B_{1''2'}\tens A_{2'''3''}\tens B_{2^43'} ,
\end{multline}
\begin{equation}\label{e66b}
\e^{A\barten B} : A^{14}\tens B^{23} \xra{A\tens\e} A^{13}\tens I^2 \simeq
A^{12} \xra\e I^1 .
\end{equation}
The tensor product $A\barten B$ is associative and turn $\Coalgsq(\vhat)$
into a monoidal category.
\end{proposition}

Graphical notation and explanation of the comultiplication~\eqref{e66a} and
the counit~\eqref{e66b} is the following
\[
\unitlength 0.70mm
\linethickness{0.4pt}
\begin{picture}(124.33,40)
\put(64.83,10.11){\oval(9.67,9.89)[t]}
\bezier{136}(45,9.89)(42.33,21.08)(65,20.22)
\bezier{124}(65,20.22)(85.33,20.22)(84.67,9.89)
\put(99.67,9.89){\line(0,1){20.22}}
\put(115,30.11){\line(0,-1){20.22}}
\put(30,9.89){\line(0,1){20.22}}
\put(14.67,30.11){\line(0,-1){20.22}}
\put(15,4.30){\makebox(0,0)[cb]{$1'$}}
\put(30,4.30){\makebox(0,0)[cb]{$1''$}}
\put(45,4.30){\makebox(0,0)[cb]{$2'$}}
\put(60,4.30){\makebox(0,0)[cb]{$2''$}}
\put(69.67,4.30){\makebox(0,0)[cb]{$2'''$}}
\put(84.67,4.30){\makebox(0,0)[cb]{$2^4$}}
\put(99.67,4.30){\makebox(0,0)[cb]{$3'$}}
\put(114.67,4.30){\makebox(0,0)[cb]{$3''$}}
\put(15,31.83){\makebox(0,0)[cb]{$1'$}}
\put(30,31.83){\makebox(0,0)[cb]{$1''$}}
\put(99.67,31.83){\makebox(0,0)[cb]{$3'$}}
\put(114.67,31.83){\makebox(0,0)[cb]{$3''$}}
\put(37.67,12.04){\makebox(0,0)[cc]{$B$}}
\put(92.33,12.04){\makebox(0,0)[cc]{$B$}}
\put(64.67,27.10){\makebox(0,0)[cc]{$B$}}
\put(37.33,-0){\makebox(0,0)[cb]{$A$}}
\bezier{26}(14.67,0.86)(14.67,0.86)(33.33,0.86)
\bezier{26}(60,0.86)(60,0.86)(41.33,0.86)
\put(92.67,-0){\makebox(0,0)[cb]{$A$}}
\bezier{26}(115.33,0.86)(115.33,0.86)(96.67,0.86)
\bezier{26}(70,0.86)(70,0.86)(88.67,0.86)
\put(65,40){\makebox(0,0)[ct]{$A$}}
\bezier{40}(14.67,37.85)(14.67,37.85)(59.67,37.85)
\bezier{40}(115.33,37.85)(115.33,37.85)(70.33,37.85)
\put(124.33,18.06){\makebox(0,0)[cc]{,}}
\put(2.33,20.22){\makebox(0,0)[cc]{$\Delta =$}}
\end{picture}
\]
\[
\unitlength 1mm
\linethickness{0.4pt}
\begin{picture}(70.17,20)
\put(37.42,12.47){\oval(15.17,15.05)[b]}
\bezier{148}(15,12.47)(14.17,-1.08)(37.50,0)
\bezier{148}(37.50,0)(61,-1.08)(59.83,12.47)
\put(15,15.05){\makebox(0,0)[cb]{$1^1$}}
\put(30,15.05){\makebox(0,0)[cb]{$1^2$}}
\put(45,15.05){\makebox(0,0)[cb]{$1^3$}}
\put(60,15.05){\makebox(0,0)[cb]{$1^4$}}
\put(37.50,12.69){\makebox(0,0)[cc]{$B$}}
\put(37.50,20){\makebox(0,0)[ct]{$A$}}
\bezier{26}(15,18.92)(15,18.92)(34,18.92)
\put(70.17,8.82){\makebox(0,0)[cc]{.}}
\bezier{26}(60,18.92)(60,18.92)(41,18.92)
\put(10,9.89){\makebox(0,0)[rc]{$\e =$}}
\end{picture}
\]
It reflects the canonical isomorphism $j_+: Y\pti\tens X\pti \to (X\tens
Y)\pti$ which is described in \secref{S12.8} on page~\pageref{p12.8}.

\begin{remark}\label{R69}
Let $A = X\odot X\pti$, $B = Y\odot Y\pti$. Then $A\barten B = (X\tens
Y)\odot(Y\pti\tens X\pti) \xra{1\odot j_+} (X\tens Y)\odot(X\tens Y)\pti$
is a squared coalgebra isomorphism.
\end{remark}

\begin{proposition}\label{P69}
Let $M\in {}^A\vhat$, $N\in {}^B\vhat$. Then $M\tens N$ has the structure
of an $A\barten B$-comodule
\begin{multline}\label{e69a}
\delta^{M\tens N} = \bigl(  M_{1'}\tens N_{1''}\odot I_2
\xra{M\tens\delta^N} M_{1'}\tens B_{1''2'}\tens N_{2''} \\
\simeq M_{1'}\tens B_{1''2'}\tens I_{2''}\tens N_{2'''} \xra{\delta^M\tens
B\tens N} A_{1'2''}\tens B_{1''2'}\tens M_{2'''}\tens N_{2^4} \bigr) .
\end{multline}
For $L\in {}^C\vhat$, both ways to construct an $A\barten B\barten
C$-comodule structure on $M\tens N\tens L$ coincide.
\end{proposition}

Graphical notation for the coaction $\delta^{M\tens N}$ is
\[
\unitlength 0.70mm
\linethickness{0.4pt}
\begin{picture}(97,40.12)
\put(64.83,10.11){\oval(9.67,9.89)[t]}
\bezier{136}(45,9.89)(42.33,21.08)(65,20.22)
\bezier{124}(65,20.22)(85.33,20.22)(84.67,9.89)
\put(30,9.89){\line(0,1){20.22}}
\put(14.67,30.11){\line(0,-1){20.22}}
\put(15,4.30){\makebox(0,0)[cb]{$1'$}}
\put(30,4.30){\makebox(0,0)[cb]{$1''$}}
\put(45,4.30){\makebox(0,0)[cb]{$2'$}}
\put(60,4.30){\makebox(0,0)[cb]{$2''$}}
\put(69.67,4.30){\makebox(0,0)[cb]{$2'''$}}
\put(84.67,4.30){\makebox(0,0)[cb]{$2^4$}}
\put(15,31.83){\makebox(0,0)[cb]{$1'$}}
\put(30,31.83){\makebox(0,0)[cb]{$1''$}}
\put(37.67,12.04){\makebox(0,0)[cc]{$B$}}
\put(37.33,-0){\makebox(0,0)[cb]{$A$}}
\bezier{26}(14.67,0.86)(14.67,0.86)(33.33,0.86)
\bezier{26}(60,0.86)(60,0.86)(41.33,0.86)
\put(2.33,20.22){\makebox(0,0)[cc]{$\delta^{M\tens N} =$}}
\put(15.02,40.12){\makebox(0,0)[ct]{$M$}}
\put(29.98,40.12){\makebox(0,0)[ct]{$N$}}
\put(69.86,-0.03){\makebox(0,0)[cb]{$M$}}
\put(85.02,-0.03){\makebox(0,0)[cb]{$N$}}
\put(97,19.35){\makebox(0,0)[cc]{.}}
\end{picture}
\]

\begin{remark}\label{R71.1}
If $A=M\odot M\pti$, $B=N\odot N\pti$, the coaction $\delta^{M\tens N}$ of
$A\barten B$ on $M\tens N$ is mapped by $1\odot j_+$ to the canonical
coaction of $(M\tens N)\odot(M\tens N)\pti$ (see \remref{R69}).
\end{remark}

\subsection{Bicoalgebras}
\begin{definition}
A {\em squared bicoalgebra} $B = (B,\Delta,\e,m,\eta)$ in $\vhat$ is
a squared coalgebra $(B,\Delta,\e) \in \Coalgsq(\vhat)$ equipped with
an algebra structure $(B,m,\eta)$ in $(\vtvhat,\barten)$ (such that
the multiplication $m:B \barten B \to B \in \vtvhat$ is associative
and $\eta: I\odot I \to B\in \vtvhat$ is the unit) and such that
$m$, $\eta$ are squared coalgebra homomorphisms, that is,
\begin{equation*}
\begin{CD}
B_{1'3''}\tens B_{1''3'}\odot I_2 @>B_{1'3''}\tens\Delta_{1''23'}>>
B_{1'3''}\tens B_{1''2'}\tens B_{2''3'} \\
@Vm_{13}\odot I_2VV            @VV\simeq V \\
B_{13}\odot I_2 @. B_{1'3''}\tens B_{1''2'}\tens I_{2''}\tens B_{2'''3'} \\
@V\Delta_{123}VV    @VV\Delta_{1'2''3''}\tens B_{1''2'}\tens B_{2'''3'}V \\
B_{12'}\tens B_{2''3} @<m\tens m<\hphantom{B_{1'3''}\tens\Delta_{1''23'}}<
B_{1'2''}\tens B_{1''2'}\tens B_{2'''3''}\tens B_{2^43'}
\end{CD}
\end{equation*}
\begin{equation*}
\begin{CD}
B^{14}\tens B^{23} @>B\tens\e>> B^{13}\tens I^{2} @>\sim>> B^{12} \\
@V\ot mVV       @.      @VV\e V \\
B^{12} @.
\makebox[0pt]{$\stackrel{\e}{\hbox to 127pt {\rightarrowfill}}$}
@. I
\end{CD}
\end{equation*}
\begin{equation*}
\begin{CD}
I_1\odot I_2\odot I_3 @>I_1\odot r^{-1}_{I\:2}\odot I_3>>
I_1\odot I_{2'}\tens I_{2''}\odot I_3 \\
@V\eta_{13}\odot I_2VV          @VV\eta_{12'}\tens\eta_{2''3}V \\
B_{13}\odot I_2 @>\Delta_{123}>\hphantom{I_1\odot r^{-1}_{I\:2}\odot I_3}>
B_{12'}\tens B_{2''3}
\end{CD}
\end{equation*}
\begin{equation*}
\bigl( I\tens I \xra{\ot\eta} \ot B \xra\e I \bigr) = r_I
\end{equation*}
hold.

Morphisms of bicoalgebras are those preserving algebra and squared
coalgebra structures. The category of bicoalgebras in $\vhat$ is
denoted $\Bicoalg(\vhat)$.
\end{definition}

\begin{remark}
Definition of a squared bicoalgebra is not self-dual. One can define
dually squared bialgebras which are ordinary coalgebras and squared
algebras.
\end{remark}

For $M,N\in {}^B\vhat$ let us define a \Bd-comodule structure on
$M\tens N\in\vhat$ via \eqref{e69a}:
\begin{multline*}
M_{1'}\tens N_{1''}\odot I_2 \xra{\delta^{M\tens N}}
B_{1'2''}\tens B_{1''2'}\tens M_{2''}\tens N_{2^4} \\
\xra{m\tens M\tens N} B_{12'}\tens M_{2''}\tens N_{2'''} .
\end{multline*}
\propref{P69} shows that the tensor product of $B$-comodules is
associative. The associativity isomorphism coincides with that one of
$(\vhat,\tens)$. Therefore, the category $({}^B\vhat,\tens)$ is
monoidal with the unit object $(I,\delta_I)$,
\[ \delta_I = \bigl( I_1\odot I_2 \xra\eta B_{12} \simeq
B_{12'}\tens I_{2''} \bigr) .\]

\subsection{Monoidal reconstruction}
\begin{definition}\label{D80}
Let $(\ca,\tens)$ be a \kd-linear abelian monoidal essentially small
category. A {\em monoidal fibre functor} is a \kd-linear exact faithful
monoidal functor $(\om_\ca,\om^\ca) : \ca \to \cv$. Let the {\em category
of monoidal fibre functors} $\fm(\cv)$ have monoidal fibre functors as
objects and let morphisms from $(\om_\ca,\om^\ca) : \ca \to \cv$ to
$(\om_\cb,\om^\cb) :  \cb \to \cv$ be equivalence classes of triples
$(F,f,\phi)$, where $(F,f) :  \ca\to \cb$ is a monoidal functor and $\phi:
(\om_\ca,\om^\ca) \to (\om_\cb,\om^\cb) \circ (F,f)$ is an isomorphism of
monoidal functors. Two such triples $(F,f,\phi)$ and $(G,g,\gamma)$ are
equivalent if there is a functorial isomorphism $\zeta:  F\to G$ such that
\eqref{e35} holds. The composite of two morphisms represented by
$(F,f,\phi)$ and $(G,g,\gamma)$ is represented by $(G\circ F, Gf\circ
g_{F,F}, \gamma_F\circ\phi)$.
\end{definition}

The functor $F$ from the above triple is exact and faithful.
Forgetting the monoidal structure
we get a functor $\fm(\cv) \to \fa(\cv)$. It is faithful since $\om_\cb
f_{X,Y}$ is determined uniquely by given $F$, $\phi$ from the condition
$\phi: (\om_\ca,\om^\ca) \xra\sim (\om_\cb,\om^\cb) \circ (F,f)$, and
$\om_\cb$ is faithful.

There is a functor $\Phi: \Bicoalg(\vhat) \to \fm(\cv)$, $B\mapsto ((\cu,
\id) : {}^B\cv \to \cv)$. We want to prove that this is an equivalence.

Let $(\om_\ca, \om^\ca) : (\ca,\tens) \to (\cv,\tens)$ be a monoidal fibre
functor. Assume that $C_\ca$ is the coend~\eqref{e38a}
constructed from the functor $\om_\ca$ (classifying
coalgebra of the category $\ca$). Then, $\ca \simeq {}^{C_\ca}\cv$.

\begin{theorem}\label{T86}
\textup{(a)}
The squared coalgebra $C_\ca$ is a bicoalgebra. The monoidal
functor $\om_\ca$ admits the factorisation
\[ (\om_\ca, \om^\ca) = \bigl(\ca \xra{(F_\ca, f^\ca)} {}^{C_\ca}\cv
\xra{(\cu,\id)} \cv \bigr) ,\]
where $F_\ca$ is the equivalence from \thmref{T49} and
$f^\ca_{X,Y} = \om^\ca_{X,Y}$.

\textup{(b)}
The correspondence $(\om_\ca,\om^\ca) \mapsto C_\ca$ extends to a functor
\[ \Psi: \fm(\cv) \to \Bicoalg(\vhat) .\]
\end{theorem}

\begin{corollary}\label{C88}
The monoidal functor
$(F_\ca,f^\ca) : \ca \to \Phi(\Psi(\ca)) = {}^{C_\ca}\cv$ constructed
in \thmref{T86} together with
$\id: (\om_\ca,\om^\ca) \to (\cu,\id)\circ(F_\ca,f^\ca)$ gives an
isomorphism of functors
\[ F: \Id_{\fm(\cv)} \xra\sim \Phi\Psi .\]
\end{corollary}

\begin{proof}
We already know by \thmref{T54} that $F$ is an isomorphism of the
forgetful functor $\fm(\cv) \to \fa(\cv)$ with
$\Phi\Psi : \fm(\cv) \to \fm(\cv) \to \fa(\cv)$. \thmref{T86} implies
that this isomorphism is in $\fm(\cv)$.
\end{proof}

\begin{theorem}\label{T89}
The functors $\Phi: \Bicoalg(\vhat) \to \fm(\cv)$, $B\mapsto {}^B\cv$ and
$\Psi: \fm(\cv) \to \Bicoalg(\vhat)$, $(\om_\ca,\om^\ca) \mapsto C_\ca$
are equivalences, quasiinverse to each other.
\end{theorem}

\begin{proof}
By \thmref{T86} $\Phi$ is essentially surjective on objects. \corref{C88}
implies that $\Phi$ is full. By \thmref{T54} $\Phi$ is faithful, hence,
$\Phi$ is an equivalence. $\Psi$ is quasi-inverse to $\Phi$ by
\corref{C88}.
\end{proof}

In the case $\cv = \kk\vect$ this theorem was proved by
Schauenburg~\cite{Sch:hopf1}.

To construct explicitly the isomorphism $\Psi\Phi \to \Id_{\fm(\cv)}$
let us use the results for $\fa(\cv)$. Let $B$ be a bicoalgebra in
$\vhat$, let $\ca= {}^B\cv$, let $\om_\ca =\cu : {}^B\cv \to \cv$ and
$\om^\ca = \id$. It was shown in \propref{P56B} that there is a unique
coalgebra isomorphism $h_B: C_\ca = \Psi\Phi(B) \to B$ such that for
any $X\in {}^B\cv$ the composite $X\odot X\pti \xra{i_X} C_\ca
\xra{h_B} B$ is the canonical coalgebra morphism \eqref{e56a}.

\begin{proposition}[bicoalgebra reconstruction]
The morphism $h_B : C_\ca \to B$ is an isomorphism of bicoalgebras
giving the functorial isomor\-phism
\[ h: \Psi\Phi \to \Id_{\Bicoalg(\vhat)} .\]
\end{proposition}

\subsection{Relationship with braided bialgebras}
Let us consider the case of braided $\cv$. Then it makes sense to consider
braided (quasiclassical) bialgebras.

There is a unique functorial isomorphism $\phi: (\ot X)\tens (\ot Y) \to
\ot(X\barten Y)$, $X,Y\in \vtv$, such that for $X = A\odot B$, $Y = C\odot
D$ it equals
\[ \phi = (432)_+\sptilde : (A\tens B)\tens (C\tens D) \to (A\tens
C)\tens(D\tens B) ,\]
\[
\unitlength 0.70mm
\linethickness{0.4pt}
\begin{picture}(66.33,10)
\put(10,9.89){\line(0,-1){9.89}}
\put(25,0){\line(3,2){15}}
\put(40,0){\line(3,2){15}}
\put(25,9.89){\line(3,-1){8.83}}
\put(55,0){\line(-3,1){9}}
\put(66.33,3.66){\makebox(0,0)[cc]{.}}
\put(1.50,4.95){\makebox(0,0)[rc]{$\phi =$}}
\put(36.02,6.31){\line(3,-1){7.59}}
\end{picture}
\]
Indeed, it extends uniquely to arbitrary $X,Y \in\vtv$ via
resolutions 
\[ 0 \to M \xra{i} A\odot B \xra{\sum f_l\odot g_l} C\odot D .\]

\begin{proposition}\label{P102A}
When $\cv$ is braided, there is a monoidal functor
\[ (\ot,\phi,r_I^{-1}) : (\vtv, \barten,I\odot I) \to (\cv,\tens,I) .\]
\end{proposition}

\begin{proposition}\label{P102B}
Let $B = (B,\Delta, \e,m,\eta)$ be a bicoalgebra in $\vhat$. Denote
\begin{align*}
\bar B &= \ot B = B^{12} ,\\
\bar\Delta &= \bigl( B^{12} \simeq B^{13}\tens I^2 \xra{\Delta^{123}}
B^{12}\tens B^{34} \bigr) ,\\
\bar\eta &= \bigl(I \xra{r_I^{-1}} I\tens I \xra{\ot\eta} \ot B \bigr) ,\\
\bar m &= \bigl( B^{12}\tens B^{34} \xra\phi B^{14}\tens B^{23} \xra{\ot m}
B^{12} \bigr).
\end{align*}
Then $(\bar B,\bar\Delta,  \e,\bar m,\bar\eta)$ is a braided bialgebra in
$\vhat$.
\end{proposition}

\begin{proposition}
There exists a \kd-linear exact faithful monoidal functor
\begin{gather*}
(F,\id,\id) : {}^B\cv \to {}^{\bar B}\cv, \quad F(X,\delta_X) = (X,
\bar\delta_X), \quad F(f) = f ,\\
\bar\delta_X = \bigl( X \xra{r_X^{-1}} X\tens I \xra{\ot\delta_X} \bar
B\tens X \bigr) ,
\end{gather*}
commuting with the underlying functor.
\end{proposition}

\section{Hopf coalgebras}
In addition to standard assumptions of \secref{S13} we assume that
$X^{\vee\vee}$ is functorially isomorphic to $X\in\cv$. Let us call such
categories rigid-involutive. In particular, braided categories are
rigid-involutive. However, we don't need the braiding in $\cv$ and our
study of Hopf algebras applies to non-braided $\cv$ as well.

Let us pick a functorial isomorphism $\zeta_X: X\to X^{\vee\vee}$. Our
notions will explicitly depend on the choice of $\zeta$. Different choices
lead to isomorphic constructions.

\subsection{Opposite coalgebras}
\begin{proposition}\label{P107}
Let $p:\cp \to \cv$ be a functor and let $p\pti:
\cp^\op \to \cv$ be the composite functor, $p\pti = \vee \circ p$, that
is, $p\pti(X) = p(X)\pti$, $p\pti(f:X\to Y) = (p(f)^t : p(Y)\pti \to
p(X)\pti)$. Let $C_p$ and $C_{p\pti}$ be the corresponding squared
coalgebras. Then there is an isomorphism $z: PC_p \to C_{p\pti}
\in\vtvhat$, which satisfies the following diagram
\begin{equation}\label{d107}
\begin{CD}
pX\pti\odot pX @>pX\pti\odot\zeta>> pX\pti\odot pX^{\vee\vee} \\
@VPi_XVV  @VVi_{X\pti}V \\
PC_p @>z>\hphantom{pX\pti\odot\zeta}> C_{p\pti}
\end{CD}
\end{equation}
for each $X\in \Ob\cp$. There are unique coalgebra structures $(pX\pti\odot
pX, \Delta^\op, \e^\op)$ and $(PC_p, \Delta^\op, \e^\op)$ such that each
morphism in diagram~\eqref{d107} is a coalgebra morphism. For $M=pX$ this
is
\begin{align}
\Delta^\op &: M\pti\odot I\odot M \xra{1\odot\coev\odot1} M\pti\odot
M^{\vee\vee} \tens M\pti\odot M \notag \\
&\hspace*{3cm} \xra{1\odot\zeta^{-1}\tens1\odot1}
M\pti\odot M\tens M\pti\odot M , \label{e108a} \\
\e^\op &: M\pti\tens M \xra{1\tens\zeta} M\pti\tens M^{\vee\vee} \xra\ev I.
\label{e108b}
\end{align}
\end{proposition}

\begin{definition}
Let $C$ be a squared coalgebra in $\vhat$. The opposite coalgebra $C_\op =
(PC,\Delta^\op,\e^\op)$ is the unique coalgebra structure on $PC$ such that
$P\ii_M : M\pti\odot M \to PC$ is a squared coalgebra homomorphism for any
\Cd-comodule $M\in {}^C\cv$, where $M\pti\odot M$ is equipped with
coalgebra structure \eqref{e108a}, \eqref{e108b}.
\end{definition}

To check the existence of the opposite coalgebra notice that any coalgebra
has the form $C_p$ for some $\cp \subset {}^C\cv$ and apply \propref{P107}.

\begin{remark}\label{R109}
The duality yields an equivalence of categories
$({}^C\cv)^\op \to {}^{C_\op}\cv$,
$(M,\delta_M) \mapsto (M\pti, \delta'_{M\pti})$,
where $\delta'$ is given by
\begin{multline*}
\delta'_{M\pti} = \bigl( M\pti\odot I \xra{M\pti\odot\coev} M\pti\odot
M^{\vee\vee}\tens M\pti \xra{M\pti\odot\zeta^{-1}\tens M\pti} \\
M\pti\odot M\tens M\pti \xra{P\ii_M\tens M\pti}
(PC)_{12'}\tens M_{2''}\pti \bigr) .
\end{multline*}
This follows by \propref{P107}.
\end{remark}

Clearly, another choice of $\zeta$ gives an isomorphic coalgebra structure
in $PC$.

\begin{proposition}\label{P109.1}
Each functorial isomorphism $\zeta: X\to X^{\vee\vee}$ determines a functor
\begin{gather*}
P_\zeta : \Coalgsq(\vhat) \to \Coalgsq(\vhat), \qquad h\mapsto Ph, \\
P_\zeta(C) = C_\op = (PC,\Delta^\op,\e^\op) = (PC,\Delta^\zeta,\e^\zeta) .
\end{gather*}
All such functors are isomorphic.
\end{proposition}

\subsection{Comparison with opposite coalgebra in braided case}
If $\cv$ is braided, we have the usual notion of an opposite coalgebra. The
following proposition shows that the new notion of opposite coincides with
traditional one at the quasiclassical level.

\begin{proposition}
Let $C$ be a squared coalgebra in $\vhat$, let
$\zeta = u_1^2 : X\to X^{\vee\vee}$ and let
$\bar C_\op = (\bar C, \bar\Delta^\op)$ be the quasiclassical opposite
to $\bar C$:
\[ \bar\Delta^\op = \bigl( \bar C \xra{\bar\Delta} \bar C\tens \bar C
\xra{c} \bar C\tens \bar C \bigr) .\]
Then the isomorphism, induced by the braiding
\[ c: \bar C_\op = (C^{12},\bar\Delta^\op,\e) \to
(C^{21},\wb{\Delta^\op},\e^\op) = \wb{C_\op} \]
is a coalgebra isomorphism.
\end{proposition}

\subsection{The antipode}
Hopf coalgebras are bicoalgebras, whose categories of comodules are rigid.
However, at the moment we use another definition: Hopf coalgebras are
bicoalgebras with an antipode. The first definition will be made a result.

\begin{definition}
Let $H$ be a bicoalgebra in $\vhat$. A \emph{right antipode} in $H$
(with respect to $\zeta$) is a morphism
$\gamma' = \gamma_\zeta : H_\op\to H \in \vtvhat$ such that
\begin{align*}
&\bigl( H_{\op1''1'}\odot I_2 \xra{\Delta_{1''21'}^\op}
H_{\op1''2'}\tens H_{\op2''1'} \rEq H_{1'2''}\tens H_{\op1''2'} \notag \\
&\hspace*{4cm} \xra{H\tens\gamma'_{1''2'}}
H_{1'2''}\tens H_{1''2'} \xra{m} H_{12} \bigr) \notag \\
&= \bigl( H_{1'1''}\odot I_2 \xra{\e\odot I} I_1\odot I_2
\xra\eta H_{12} \bigr) , 
\end{align*}
\begin{align*}
&\bigl( I_1\odot H_{2''2'} \xra{\Delta_{2''12'}} H_{2''1'}\tens H_{1''2'}
\rEq H_{\op1'2''}\tens H_{1''2'} \notag \\
&\hspace*{4cm} \xra{\gamma'_{1'2''}\tens H}
H_{1'2''}\tens H_{1''2'} \xra{m} H_{12} \bigr) \notag   \\
&= \bigl( I_1\odot H_{\op2'2''} \xra{I\odot\e^\op} I_1\odot I_2
\xra\eta H_{12} \bigr) . 
\end{align*}
A \emph{left antipode} in $H$ (with respect to $\zeta$) is a morphism
${}'\gamma = {}_\zeta\gamma : H_\op\to H \in \vtvhat$ such that
\begin{gather*}
\begin{split}
&\bigl( H_{1''1'}\odot I_2 \xra{\Delta_{1''21'}} H_{1''2'}\tens H_{2''1'}
\rEq H_{\op1'2''}\tens H_{1''2'} \\
&\hspace*{4cm} \xra{{}'\gamma_{1'2''}\tens H}
H_{1'2''}\tens H_{1''2'} \xra{m} H_{12} \bigr)   \\
&= \bigl( H_{\op1'1''}\odot I_2 \xra{\e^\op\odot I_2} I_1\odot I_2
\xra\eta H_{12} \bigr) ,
\end{split}
\\
\begin{split}
&\bigl( I_1\odot H_{\op2''2'} \xra{\Delta_{2''12'}^\op}
H_{\op2''1'}\tens H_{\op1''2'} \rEq H_{1'2''}\tens H_{\op1''2'} \\
&\hspace*{4cm} \xra{H\tens{}'\gamma_{1''2'}}
H_{1'2''}\tens H_{1''2'} \xra{m} H_{12} \bigr)   \\
&= \bigl( I_1\odot H_{2'2''} \xra{I\odot\e} I_1\odot I_2
\xra\eta H_{12} \bigr) .
\end{split}
\end{gather*}
A (\emph{squared}) \emph{Hopf coalgebra} is a bicoalgebra which has
a right and a left antipode.
\end{definition}

Graphical expression of these equations is the following. Here $X$ is an
\Hd-comodule and $\ii_X : X\odot X\pti \to H$ is implicit.
\begin{gather*}
\unitlength 0.50mm
\linethickness{0.4pt}
\begin{picture}(140,89.90)
\put(0,9.68){\framebox(70,10.32)[cc]{$m$}}
\put(20,24.95){\framebox(30,10.11)[cc]{$\gamma'$}}
\put(45,24.95){\line(0,-1){4.95}}
\put(25,24.95){\line(0,-1){4.95}}
\put(45,35.05){\line(0,1){4.95}}
\put(55,50){\oval(20,20)[t]}
\put(65,49.89){\line(0,-1){29.89}}
\put(55,63.01){\makebox(0,0)[cb]{2}}
\put(68,55.05){\makebox(0,0)[lc]{$X\pti$}}
\put(44.67,55.05){\makebox(0,0)[rc]{$X^{\vee\vee}$}}
\put(25,64.95){\line(0,-1){29.89}}
\put(5,64.95){\line(0,-1){44.95}}
\put(15,9.68){\line(0,-1){4.73}}
\put(55,9.68){\line(0,-1){4.73}}
\put(60.17,3.01){\makebox(0,0)[cc]{2}}
\put(10,3.01){\makebox(0,0)[cc]{1}}
\put(35,-0){\makebox(0,0)[cb]{$H$}}
\put(2,60){\makebox(0,0)[rc]{$1'$}}
\put(22,60){\makebox(0,0)[rc]{$1''$}}
\put(25,67.96){\makebox(0,0)[cb]{$X\pti$}}
\put(5,67.96){\makebox(0,0)[cb]{$X$}}
\put(110,9.89){\framebox(30,10.32)[cc]{$\eta$}}
\put(125,64.95){\oval(20,49.89)[b]}
\put(90,31.83){\makebox(0,0)[cc]{=}}
\put(135.02,9.87){\line(0,-1){4.94}}
\put(115.02,9.87){\line(0,-1){4.94}}
\put(111.89,4.94){\makebox(0,0)[rc]{1}}
\put(138.15,4.94){\makebox(0,0)[lc]{2}}
\put(124.76,0){\makebox(0,0)[cb]{$H$}}
\put(135,67.96){\makebox(0,0)[cb]{$X\pti$}}
\put(115,67.96){\makebox(0,0)[cb]{$X$}}
\put(38,40.05){\framebox(14.10,9.82)[cc]{$\zeta^{-1}$}}
\end{picture}
\label{f112i} \\[5mm]
\unitlength 0.50mm
\linethickness{0.4pt}
\begin{picture}(140.05,53)
\put(0,9.68){\framebox(70,10.32)[cc]{$m$}}
\put(15,9.68){\line(0,-1){4.73}}
\put(55,9.68){\line(0,-1){4.73}}
\put(60.17,3.01){\makebox(0,0)[cc]{2}}
\put(10,3.01){\makebox(0,0)[cc]{1}}
\put(35,-0){\makebox(0,0)[cb]{$H$}}
\put(110,9.89){\framebox(30,10.32)[cc]{$\eta$}}
\put(135.02,9.87){\line(0,-1){4.94}}
\put(115.02,9.87){\line(0,-1){4.94}}
\put(111.89,4.94){\makebox(0,0)[rc]{1}}
\put(138.15,4.94){\makebox(0,0)[lc]{2}}
\put(124.76,0){\makebox(0,0)[cb]{$H$}}
\put(0.22,25.01){\framebox(9.78,10.10)[cc]{$\gamma'$}}
\put(59.97,25.01){\framebox(9.96,10.10)[cc]{}}
\bezier{30}(59.97,30.06)(35.25,30.06)(10,30.06)
\put(14.98,34.99){\oval(19.92,19.96)[t]}
\put(5.02,25.01){\line(0,-1){5.05}}
\put(24.94,35.11){\line(0,-1){15.14}}
\put(64.95,25.01){\line(0,-1){5.05}}
\put(45.03,19.96){\line(0,1){30.06}}
\put(64.95,35.11){\line(0,1){14.91}}
\put(2,41.99){\makebox(0,0)[rc]{$X\pti$}}
\put(27.96,41.99){\makebox(0,0)[lc]{$X$}}
\put(14.98,47.96){\makebox(0,0)[cb]{1}}
\put(48.01,47.96){\makebox(0,0)[lc]{$2'$}}
\put(67.97,47.96){\makebox(0,0)[lc]{$2''$}}
\put(64.95,53){\makebox(0,0)[cb]{$X$}}
\put(45.03,53){\makebox(0,0)[cb]{$X\pti$}}
\put(130.07,34.95){\framebox(9.98,10.01)[cc]{$\zeta$}}
\put(124.99,34.95){\oval(20.14,20.03)[b]}
\put(135.06,49.97){\line(0,-1){5.01}}
\put(114.91,34.95){\line(0,1){15.02}}
\put(135.28,53){\makebox(0,0)[cb]{$X$}}
\put(115.37,53){\makebox(0,0)[cb]{$X\pti$}}
\put(89.67,27.96){\makebox(0,0)[cc]{=}}
\end{picture}
\label{f112ii}
\end{gather*}
\begin{gather*}
\unitlength 0.50mm
\linethickness{0.4pt}
\begin{picture}(140.05,53)
\put(0,9.68){\framebox(70,10.32)[cc]{$m$}}
\put(15,9.68){\line(0,-1){4.73}}
\put(55,9.68){\line(0,-1){4.73}}
\put(60.17,3.01){\makebox(0,0)[cc]{2}}
\put(10,3.01){\makebox(0,0)[cc]{1}}
\put(35,-0){\makebox(0,0)[cb]{$H$}}
\put(110,9.89){\framebox(30,10.32)[cc]{$\eta$}}
\put(135.02,9.87){\line(0,-1){4.94}}
\put(115.02,9.87){\line(0,-1){4.94}}
\put(111.89,4.94){\makebox(0,0)[rc]{1}}
\put(138.15,4.94){\makebox(0,0)[lc]{2}}
\put(124.76,0){\makebox(0,0)[cb]{$H$}}
\put(0.22,25.01){\framebox(9.78,10.10)[cc]{${}'\gamma$}}
\put(59.97,25.01){\framebox(9.96,10.10)[cc]{}}
\bezier{30}(59.97,30.06)(35.25,30.06)(10,30.06)
\put(54.97,34.99){\oval(19.92,19.96)[t]}
\put(5.02,25.01){\line(0,-1){5.05}}
\put(64.95,25.01){\line(0,-1){5.05}}
\put(44.66,41.66){\makebox(0,0)[rc]{$X\pti$}}
\put(67.95,41.99){\makebox(0,0)[lc]{$X$}}
\put(54.97,47.96){\makebox(0,0)[cb]{2}}
\put(2.02,47.96){\makebox(0,0)[rc]{$1'$}}
\put(21.98,47.96){\makebox(0,0)[rc]{$1''$}}
\put(24.96,53){\makebox(0,0)[cb]{$X$}}
\put(5.04,53){\makebox(0,0)[cb]{$X\pti$}}
\put(130.07,34.95){\framebox(9.98,10.01)[cc]{$\zeta$}}
\put(124.99,34.95){\oval(20.14,20.03)[b]}
\put(135.06,49.97){\line(0,-1){5.01}}
\put(114.91,34.95){\line(0,1){15.02}}
\put(135.28,53){\makebox(0,0)[cb]{$X$}}
\put(115.37,53){\makebox(0,0)[cb]{$X\pti$}}
\put(89.67,27.96){\makebox(0,0)[cc]{=}}
\put(45.07,34.98){\line(0,-1){15.08}}
\put(5.08,50.06){\line(0,-1){15.08}}
\put(24.99,50.06){\line(0,-1){30.15}}
\end{picture}
\label{f113iii} \\
\unitlength 0.50mm
\linethickness{0.4pt}
\begin{picture}(140,69.17)
\put(0,9.68){\framebox(70,10.32)[cc]{$m$}}
\put(20,24.95){\framebox(30,10.11)[cc]{${}'\gamma$}}
\put(45,24.95){\line(0,-1){4.95}}
\put(25,24.95){\line(0,-1){4.95}}
\put(15,44.84){\oval(20,20)[t]}
\put(15,57.85){\makebox(0,0)[cb]{1}}
\put(28,49.89){\makebox(0,0)[lc]{$X\pti$}}
\put(2,49.89){\makebox(0,0)[rc]{$X^{\vee\vee}$}}
\put(15,9.68){\line(0,-1){4.73}}
\put(55,9.68){\line(0,-1){4.73}}
\put(60.17,3.01){\makebox(0,0)[cc]{2}}
\put(10,3.01){\makebox(0,0)[cc]{1}}
\put(35,-0){\makebox(0,0)[cb]{$H$}}
\put(47.88,54.97){\makebox(0,0)[lc]{$2'$}}
\put(67.88,54.97){\makebox(0,0)[lc]{$2''$}}
\put(64.88,62.93){\makebox(0,0)[cb]{$X\pti$}}
\put(44.88,62.93){\makebox(0,0)[cb]{$X$}}
\put(110,9.89){\framebox(30,10.32)[cc]{$\eta$}}
\put(125,60.22){\oval(20,49.89)[b]}
\put(90,31.83){\makebox(0,0)[cc]{=}}
\put(135.02,9.87){\line(0,-1){4.94}}
\put(115.02,9.87){\line(0,-1){4.94}}
\put(111.89,4.94){\makebox(0,0)[rc]{1}}
\put(138.15,4.94){\makebox(0,0)[lc]{2}}
\put(124.76,0){\makebox(0,0)[cb]{$H$}}
\put(135,63.23){\makebox(0,0)[cb]{$X\pti$}}
\put(115,63.23){\makebox(0,0)[cb]{$X$}}
\put(5,34.84){\line(0,-1){14.62}}
\put(25.05,45.05){\line(0,-1){10.05}}
\put(45.08,59.89){\line(0,-1){24.89}}
\put(64.93,59.89){\line(0,-1){39.97}}
\put(-2.06,35.14){\framebox(14.10,9.82)[cc]{$\zeta^{-1}$}}
\end{picture}
\label{f113iv}
\end{gather*}
The dual notion is called squared Hopf algebra.

\begin{proposition}\label{P113}
Let $(H,\gamma', {}'\gamma)$ be a Hopf coalgebra. Then
$\gamma', {}'\gamma :H_\op \to H$ are homomorphisms of squared coalgebras.
\end{proposition}

\begin{theorem}\label{T118}
The category ${}^H\cv$ of comodules over a Hopf coalgebra $H$ is rigid.
For the right dual of $X$ one can take
\[ \bigl( X\pti, \delta_{X\pti}: X_1\pti\odot I_2 \xra{\delta'_{X\pti}}
H_{\op12'}\tens X\pti_{2''} \xra{\gamma'_{12'}\tens X\pti_{2''}}
H_{12'}\tens X\pti_{2''} \bigr) ,\]
for the left dual --
\[ (\lpti X, \delta_{\lpti X} : \lpti X_1\odot I_2
\xra{{}'\delta_{\lpti X}} H_{\op12'}\tens\lpti X_{2''}
\xra{{}'\gamma_{12'}\tens\lpti X_{2''}} H_{12'}\tens\lpti X_{2''} \bigr).
\]
The evaluation and coevaluation morphisms are the same as in $\cv$.
\end{theorem}

The formulae 
\begin{multline*}
\delta_{X\pti} = \bigl( X\pti\odot I \xra{X\pti\odot\coev}
X\pti\odot X^{\vee\vee}\tens X\pti \xra{X\pti\tens\zeta^{-1}\tens X\pti}
X_1\pti\odot X_{2'}\tens X_{2''}\pti \\
\xra{\ii_{2'1}\tens X\pti} H_{\op12'}\tens X_{2''}\pti
\xra{\gamma'\tens X\pti} H_{12'}\tens X_{2''}\pti \bigr),
\end{multline*}
\begin{multline*}
\delta_{\lpti X} = \bigl( \lpti X\odot I \xra{\lpti X\odot\coev}
\lpti X\odot X\tens\lpti X \xra{{}^t\zeta^{-1}\odot X\tens\lpti X}
X_1\pti\odot X_{2'}\tens\lpti X_{2''} \\
\xra{\ii_{2'1}\tens\lpti X} H_{\op12'}\tens\lpti X_{2''}
\xra{{}'\gamma\tens\lpti X} H_{12'}\tens\lpti X_{2''} \bigr).
\end{multline*}
together with
definition~\eqref{e56a} of $\ii_X$ imply
\begin{gather}
\begin{CD}
X_1\pti\odot X_2 @>X\pti\odot\zeta>> X_1\pti\odot X_2^{\vee\vee} \\
@VP\ii_XVV  @VV\ii_{X\pti}V \\
H_{\op12} @>\gamma'>\hphantom{X\pti\odot\zeta}> H_{12}
\end{CD}
\quad , \label{d124a} \\
\begin{CD}
X_1\pti\odot X_2 @>{}^t\zeta\odot X>> \lpti X_1\odot X_2 \\
@VP\ii_XVV  @VV\ii_{\lpti X}V \\
H_{\op12} @>{}'\gamma>\hphantom{{}^t\zeta\odot X}> H_{12}
\end{CD}
\quad . \label{d124b}
\end{gather}

\begin{corollary}
The right and the left antipodes $\gamma', {}'\gamma : H_\op \to H$ of a
Hopf coalgebra are unique and invertible.
\end{corollary}

Uniqueness follows from the fundamental theorem on coalgebras
(\corref{C58}) and diagrams \eqref{d124a}, \eqref{d124b}. Invertibility
follows from the fact that $X\mapsto X\pti$, $X\mapsto\lpti X$ are
equivalences. If only right antipode exists for $H$, it need not be
invertible.

\begin{remark}\label{R124.1}
The property of being a Hopf
coalgebra does not depend on the choice of $\zeta$.
\end{remark}

\begin{proposition}\label{P125}
Let $H$ be a bicoalgebra in $\vhat$.
\begin{enumerate}
\renewcommand{\labelenumi}{(\alph{enumi})}
\item If $H$ has a right antipode $\gamma_\zeta: H_\op \to H$ with respect
to $\zeta$ and it is invertible, then ${}_{{}^t\zeta^{-1}}\gamma =
P\gamma_\zeta^{-1}$ is a left antipode for $H$ with respect to
${}^t\zeta^{-1}$.
\item If $H$ has a left antipode ${}_\zeta\gamma : H_\op \to H$ with
respect to $\zeta$ and it is invertible, then $\gamma_{{}^t\zeta^{-1}} =
P{}_\zeta\gamma^{-1}$ is a right antipode for $H$ with respect to
${}^t\zeta^{-1}$.
\end{enumerate}
In both cases $H$ is a Hopf coalgebra.
\end{proposition}

\subsection{Comparison with braided Hopf algebras}
We recall that if $\cv$ is braided, the tensor functor $(\ot,\phi,r_I^{-1})
: \vtv \to \cv$ from \propref{P102A} transforms squared bicoalgebras $B$
into quasiclassical bialgebras $\bar B$ (see \propref{P102B}).

\begin{proposition}
Let $H$ be a squared Hopf algebra in $\vhat$. Then $\bar H$ is a
quasiclassical Hopf algebra in $\vhat$. If $\gamma'$ is the right antipode
for $H$ with respect to $\zeta = u_1^2$, then
\[ \gamma_{\bar H} = \bigl( H^{12} \xra{c} H^{21} \xra{\ot\gamma'} H^{12}
\bigr) \]
is the antipode for $\bar H$. If ${}'\gamma$ is the left antipode for $H$
with respect to $\tilde\zeta = u_{-1}^2$, then
\[ \tilde\gamma_{\bar H} = \gamma_{\bar H}^{-1} = \bigl( H^{12}
\xra{c^{-1}} H^{21} \xra{\ot{}'\gamma} H^{12} \bigr) \]
is the skew antipode for $\bar H$.
\end{proposition}

\section{Quasitriangular Hopf coalgebras}
We want to discuss braided monoidal categories. Naturally, we assume that
$\cv$ is braided. It was shown in Propositions \ref{P102A}, \ref{P102B}
that a monoidal functor $(\ot,\phi,r_I^{-1}): (\cv\tens\cv, \barten, I\odot
I) \to (\cv,\tens,I)$ maps bicoalgebras $H$ to braided bialgebras $\bar H =
(\ot H,\bar\Delta, \e,\bar m,\bar\eta)$. The structure responsible for the
braiding in ${}^B\cv$ is the \Rd-matrix. Bicoalgebras which admit an
\Rd-matrix will be called quasitriangular. We are going to establish
correspondence between \kd-linear exact monoidal functors $\omega: \cc \to
\cv$ which do not preserve braiding and quasitriangular bicoalgebras.
Quasitriangular Hopf coalgebra is a natural implementation of the idea of a
quantum braided group \cite{Hla:qbg}. Although there are two \Rd-matrices
-- $R_+$ and $R_-$, they are linearly expressed one through another with
the help of braiding in $\cv$.

\subsection{$R$-matrices in Hopf coalgebras}
\begin{definition}\label{D159}
A \emph{quasitriangular Hopf coalgebra} in $\vhat$ is a Hopf coalgebra $H$
in $\vhat$ equipped with bilinear forms $R_+: \bar H\tens\bar H \to I$,
$R_-: \bar H\tens\bar H \to I \in \Mor\vhat$ called the \Rd-matrices
such that
\begin{equation}\label{e160a}
R_+ = \bigl( \bar H\tens\bar H \xra\Omega \bar H\tens\bar H
\xra{R_-} I \bigr) ,
\end{equation}
where $\Omega$ denotes the morphism $1^1\tens(c^{32}\circ c^{23})\tens1^4:
H^{12}\tens H^{34} \to H^{12}\tens H^{34}$,
\begin{equation*}\label{e160b}
\e = \bigl( \bar H \simeq I\tens\bar H \xra{\bar\eta\tens\bar H}
\bar H\tens\bar H \xra{R_+} I \bigr) ,
\end{equation*}
\begin{equation*}\label{e160c}
\e = \bigl( \bar H \simeq \bar H\tens I \xra{\bar H\tens\bar\eta}
\bar H\tens\bar H \xra{R_+} I \bigr) ,
\end{equation*}
\begin{subequations}
\begin{equation}\label{e160d}
\begin{split}
&\bigl( \bar H\tens\bar H\tens\bar H \xra{\bar m\tens\bar H}
\bar H\tens\bar H \xra{R_+} I \bigr) \\
&= \bigl( \bar H\tens\bar H\tens\bar H
\xra{\bar H\tens\bar H\tens\bar\Delta}
\bar H\tens\bar H\tens\bar H\tens\bar H \\
&\hspace*{1cm} \xra{\bar H\tens R_+\tens\bar H}
\bar H\tens I\tens\bar H \simeq \bar H\tens\bar H \xra{R_+} I \bigr) ,
\end{split}
\end{equation}
\begin{equation}\label{e160e}
\begin{split}
&\bigl( \bar H\tens\bar H\tens\bar H \xra{\bar H\tens\bar m}
\bar H\tens\bar H \xra{R_-} I \bigr) \\
&= \bigl( \bar H\tens\bar H\tens\bar H
\xra{\bar\Delta\tens\bar H\tens\bar H}
\bar H\tens\bar H\tens\bar H\tens\bar H \\
&\hspace*{1cm} \xra{\bar H\tens c^{-1}\tens\bar H}
\bar H\tens\bar H\tens\bar H\tens\bar H \xra{R_-\tens R_-}
I\tens I \simeq I \bigr) ,
\end{split}
\end{equation}
\end{subequations}
\begin{equation}\label{e160f}
\begin{split}
&\bigl( H_{1'2'}\tens H_{1''2''} \simeq
H_{1'2'}\tens H_{1''2''}\tens I_{1'''} \xra{\Delta_{1'1'''2'}\tens H}
H_{1'1'''}\tens H_{1''2''}\tens H_{1^42'} \\
&\hspace*{1cm} \xra{c_{1''1'''}^{-1}\tens l_H^{-1}}
H_{1'1''}\tens H_{1'''2''}\tens I_{1^4}\tens H_{1^52'}
\xra{H\tens\Delta_{1'''1^42''}\tens H} \\
&\hspace*{1cm} H_{1'1''}\tens H_{1'''1^4}\tens H_{1^52''}\tens H_{1^62'}
\xra{R_+\tens m} I_{1'}\tens H_{1''2} \simeq H_{12} \bigr) \\
&= \bigl( H_{1'2'}\tens H_{1''2''} \simeq
H_{1'2''}\tens I_{2'}\tens H_{1''2'''} \xra{H\tens\Delta_{1''2'2'''}}
H_{1'2'''}\tens H_{1''2'}\tens H_{2''2^4} \\
&\hspace*{1cm} \xra{r_H^{-1}\tens c_{2''2'''}^{-1}}
H_{1'2'''}\tens H_{1''2'}\tens I_{2''}\tens H_{2^42^5}
\xra{\Delta_{1'2''2'''}\tens H\tens H} \\
&\hspace*{1cm} H_{1'2''}\tens H_{1''2'}\tens H_{2'''2^4}\tens H_{2^52^6}
\xra{m\tens R_+} H_{12'}\tens I_{2''} \simeq H_{12} \bigr) .
\end{split}
\end{equation}
\end{definition}

It is easier to look at these  conditions in graphical form.
Equation~\eqref{e160d} becomes
\begin{equation*}\label{e161a}
\unitlength 0.5mm
\linethickness{0.4pt}
\begin{picture}(160,54)
\put(0,25){\framebox(30,10)[cc]{$\bar m$}}
\put(25,45){\line(0,-1){10}}
\put(5,35){\line(0,1){10}}
\put(5,48){\makebox(0,0)[cb]{$\bar H$}}
\put(45,45){\line(0,-1){30}}
\put(15,25){\line(0,-1){10}}
\put(10,5){\framebox(40,10)[cc]{$R_+$}}
\put(25,48){\makebox(0,0)[cb]{$\bar H$}}
\put(45,48){\makebox(0,0)[cb]{$\bar H$}}
\put(90,0){\framebox(70,10)[cc]{$R_+$}}
\put(110,15){\framebox(30,10)[cc]{$R_+$}}
\put(130,30){\framebox(30,10)[cc]{$\bar\Delta$}}
\put(155,30){\line(0,-1){20}}
\put(135,25){\line(0,1){5}}
\put(145,40){\line(0,1){5}}
\put(145,48){\makebox(0,0)[cb]{$\bar H$}}
\put(115,45){\line(0,-1){20}}
\put(95,10){\line(0,1){35}}
\put(95,48){\makebox(0,0)[cb]{$\bar H$}}
\put(115,48){\makebox(0,0)[cb]{$\bar H$}}
\put(70,27){\makebox(0,0)[cc]{=}}
\end{picture}
\end{equation*}
\eqref{e160e} becomes
\begin{equation*}\label{e161b}
\unitlength 0.5mm
\linethickness{0.4pt}
\begin{picture}(169.67,44)
\put(0,0){\framebox(30,10)[cc]{$R_-$}}
\put(50,0){\framebox(30,10)[cc]{$R_-$}}
\put(55,10){\line(-3,1){30}}
\put(25,10){\line(3,1){12}}
\put(43,16){\line(3,1){12}}
\put(55,20){\line(0,1){15}}
\put(75,35){\line(0,-1){25}}
\put(75,38){\makebox(0,0)[cb]{$\bar H$}}
\put(55,38){\makebox(0,0)[cb]{$\bar H$}}
\put(0,20){\framebox(30,10)[cc]{$\bar\Delta$}}
\put(15,30){\line(0,1){5}}
\put(15,38){\makebox(0,0)[cb]{$\bar H$}}
\put(5,20){\line(0,-1){10}}
\put(100,21){\makebox(0,0)[cc]{=}}
\put(125,35){\line(0,-1){25}}
\put(125,38){\makebox(0,0)[cb]{$\bar H$}}
\put(145,35){\line(0,-1){5}}
\put(145,38){\makebox(0,0)[cb]{$\bar H$}}
\put(165,35){\line(0,-1){5}}
\put(165,38){\makebox(0,0)[cb]{$\bar H$}}
\put(140,20){\framebox(29.67,10)[cc]{$\bar m$}}
\put(155,20){\line(0,-1){10}}
\put(120,0){\framebox(40,10)[cc]{$R_-$}}
\end{picture}
\end{equation*}
and \eqref{e160f} becomes
\begin{equation*}\label{e161c}
\unitlength 0.37mm
\linethickness{0.4pt}
\begin{picture}(320,51)
\put(75,23){\oval(20,20)[t]}
\put(84,32){\makebox(0,0)[lb]{$Y$}}
\put(68,32){\makebox(0,0)[rb]{$Y\pti$}}
\put(75,23){\oval(60,40)[rt]}
\put(103,41){\makebox(0,0)[lb]{$X$}}
\put(125,23){\line(0,1){25}}
\put(145,48){\line(0,-1){25}}
\put(145,51){\makebox(0,0)[cb]{$Y\pti$}}
\put(125,51){\makebox(0,0)[cb]{$X\pti$}}
\put(135,10){\line(0,-1){5}}
\put(95,10){\line(0,-1){5}}
\put(115,0){\makebox(0,0)[cb]{$H$}}
\put(45,23){\line(0,1){25}}
\put(45,51){\makebox(0,0)[cb]{$Y$}}
\put(5,48){\line(0,-1){25}}
\put(5,51){\makebox(0,0)[cb]{$X$}}
\bezier{116}(42.67,38.67)(24.33,32.33)(25,23)
\bezier{108}(48,40.33)(61.33,43)(75,43)
\put(29.67,34.67){\makebox(0,0)[rb]{$X\pti$}}
\put(159.33,30){\makebox(0,0)[cc]{=}}
\put(0,10){\framebox(70,13)[cc]{$R_+$}}
\put(80,10){\framebox(70,13)[cc]{$m$}}
\put(225,10){\line(0,-1){5}}
\put(185,10){\line(0,-1){5}}
\put(205,0){\makebox(0,0)[cb]{$H$}}
\put(170,10){\framebox(70,13)[cc]{$m$}}
\put(250,10){\framebox(70,13)[cc]{$R_+$}}
\put(175,48){\line(0,-1){25}}
\put(195,23){\line(0,1){25}}
\put(194.67,51){\makebox(0,0)[cb]{$Y$}}
\put(175,51){\makebox(0,0)[cb]{$X$}}
\put(245,23){\oval(20,20)[t]}
\put(254,32){\makebox(0,0)[lb]{$X$}}
\put(238,32){\makebox(0,0)[rb]{$X\pti$}}
\put(255,23){\oval(80,40)[t]}
\put(291.67,39.33){\makebox(0,0)[lb]{$Y$}}
\put(219,39.33){\makebox(0,0)[rb]{$Y\pti$}}
\put(315,48){\line(0,-1){25}}
\put(315,51){\makebox(0,0)[cb]{$Y\pti$}}
\put(275,48){\line(0,-1){3.33}}
\put(275,40){\line(0,-1){17}}
\put(275,51){\makebox(0,0)[cb]{$X\pti$}}
\end{picture}
\end{equation*}
With a certain effort one can see that these properties are the dual ones
to the equations for \Rd-matrix written by Drinfeld~\cite{Dri:ICM}.

\begin{theorem}\label{T161}
Let $(H,R_+,R_-)$ be a quasitriangular Hopf coalgebra. Then the categories
of comodules ${}^H\cv$ and ${}^H\vhat$ are braided and the braiding
$R_{X,Y}: X\tens Y \to Y\tens X$ is given by
\[
\begin{split}
R_{X,Y} &= \bigl( X\tens Y \xra{\bar\delta_X\tens\bar\delta_Y}
\bar H\tens X\tens\bar H\tens Y \xra{(432)_+\sptilde}
\bar H\tens \bar H\tens Y\tens X \\
&\hspace*{5cm} \xra{R_+\tens Y\tens X} I\tens Y\tens X \simeq
Y\tens X \bigr) \\
&= \bigl( X\tens Y \xra{\bar\delta_X\tens\bar\delta_Y}
\bar H\tens X\tens\bar H\tens Y \xra{(432)_-\sptilde}
\bar H\tens \bar H\tens Y\tens X \\
&\hspace*{5cm} \xra{R_-\tens Y\tens X} I\tens Y\tens X \simeq
Y\tens X \bigr) ,
\end{split}
\]
or
\begin{equation}\label{f162b}
\raisebox{-12mm}{
\unitlength 0.40mm
\linethickness{0.4pt}
\begin{picture}(220,60)
\put(30,10){\framebox(30,10)[cc]{$R_+$}}
\put(30,40){\framebox(30,12)[cc]{$\bar\delta_X$}}
\put(35,40){\line(0,-1){20}}
\put(32,23){\makebox(0,0)[rb]{$\bar H$}}
\put(55,20){\line(3,2){30}}
\put(80,40){\framebox(30,12)[cc]{$\bar\delta_Y$}}
\put(105,40){\line(-5,-6){25}}
\put(55,40){\line(5,-3){13.33}}
\put(105,10){\line(-5,3){14.33}}
\put(86,21.67){\line(-5,3){12}}
\put(80,0){\makebox(0,0)[cb]{$Y$}}
\put(105,0){\makebox(0,0)[cb]{$X$}}
\put(55.67,23){\makebox(0,0)[rb]{$\bar H$}}
\put(45,52){\line(0,1){5}}
\put(45,60){\makebox(0,0)[cb]{$X$}}
\put(95,52){\line(0,1){5}}
\put(95,60){\makebox(0,0)[cb]{$Y$}}
\put(120,32){\makebox(0,0)[cc]{$=$}}
\put(15,32){\makebox(0,0)[rc]{$R_{X,Y} =$}}
\put(140,10){\framebox(30,10)[cc]{$R_+$}}
\put(140,40){\framebox(30,12)[cc]{$\bar\delta_X$}}
\put(145,40){\line(0,-1){20}}
\put(142,23){\makebox(0,0)[rb]{$\bar H$}}
\put(190,40){\framebox(30,12)[cc]{$\bar\delta_Y$}}
\put(185,0){\makebox(0,0)[cb]{$Y$}}
\put(215,0){\makebox(0,0)[cb]{$X$}}
\put(165.67,23){\makebox(0,0)[rb]{$\bar H$}}
\put(155,52){\line(0,1){5}}
\put(155,60){\makebox(0,0)[cb]{$X$}}
\put(205,52){\line(0,1){5}}
\put(205,60){\makebox(0,0)[cb]{$Y$}}
\put(165,40){\line(5,-3){50}}
\put(195,40){\line(-3,-2){12}}
\put(165,20){\line(3,2){12}}
\put(215,40){\line(-1,-1){16}}
\put(185,10){\line(1,1){10}}
\end{picture}
} \quad.
\end{equation}
\end{theorem}

\begin{remark}
One can define quasitriangular bicoalgebras in the same way as in
\defref{D159}, adding one condition -- invertibility of \eqref{f162b}.
\end{remark}

\begin{remark}
Equation \eqref{e160f} can be written in 4 equivalent forms. By
\eqref{e160a} one can replace $c^{-1}$, $R_+$ by $c$, $R_-$ in the left
and/or, independently, in the right hand side of \eqref{e160f}.
\end{remark}

\begin{proposition}\label{P166}
Let $H$ be a bicoalgebra (e.g. a Hopf coalgebra). If ${}^H\cv$ has a
braiding $R_{X,Y}: X\tens Y \to Y\tens X$, then $(H, R_+, R_-)$ is
quasitriangular, where $R_+$, $R_-$ are determined by
\begin{equation*}\label{e166}
\begin{split}
&\bigl( X\tens X\pti\tens Y\tens Y\pti \xra{\ii_X\tens\ii_Y} \bar
H\tens\bar H \xra{R_\pm} I \bigr) \\
&= \bigl( X\tens X\pti\tens Y\tens Y\pti \xra{X\tens c^{\pm1}\tens Y\pti}
X\tens Y\tens X\pti\tens Y\pti \xra{R_{X,Y}\tens X\pti\tens Y\pti} \\
&\hspace*{3mm} Y\tens X\tens X\pti\tens Y\pti \xra{Y\tens\ev\tens Y\pti}
Y\tens I\tens Y\pti \simeq Y\tens Y\pti \xra\ev I \bigr) .
\end{split}
\end{equation*}
\end{proposition}

\subsection{Braiding for comodules over a braided Hopf algebra}
It seems that there is no gadget, which would
make the \emph{whole} category of comodules over a braided Hopf algebra
into a braided category. That is why the framework of braided Hopf algebras
seems insufficient for the ideas like quantum braided groups.
However, braided Hopf algebras, which came from quasitriangular squared Hopf
coalgebras, make an exception.

Let $H$ be a quasitriangular Hopf coalgebra, and let $\bar H = \ot H$. The
category ${}^{\bar H}\cv$ of comodules over the braided Hopf algebra $\bar
H$ is equivalent to the tensor product ${}^H\cv\tens\cv$ of two braided
rigid monoidal categories by \thmref{T58.2}. This allows to introduce a
braiding in ${}^{\bar H}\cv$.

Indeed, the tensor product
\begin{multline*}
\unten: ({}^H\cv\tens\cv) \times ({}^H\cv\tens\cv) \xra\odot
{}^H\cv\tens\cv\tens{}^H\cv\tens\cv \\
\xra{1\tens P\tens1} {}^H\cv\tens{}^H\cv\tens\cv\tens\cv \xra{\ot\tens\ot}
{}^H\cv\tens\cv ,
\end{multline*}
\begin{alignat*}2
(A\unten B)_{12} &= A_{1'2'}\tens B_{1''2''} &\qquad& \text{for } A,B \in
{}^H\cv\tens\cv , \\
(M\odot X)\unten(N\odot Y) &= (M\tens N)\odot(X\tens Y) &\qquad&
\text{for } M,N\in {}^H\cv,\ X,Y\in \cv ,
\end{alignat*}
with the corresponding associativity and unit constraints makes
$({}^H\cv\tens\cv, \unten, I\odot I)$ into a monoidal category similarly to
\secref{S15.4}.

\begin{proposition}
The functorial isomorphism $\psi: A^{12}\tens B^{34} \to A^{13}\tens
B^{24}$, $A,B\in {}^H\cv\tens\cv$, determined by
\begin{equation*}\label{e166.2a}
\psi = 1\tens c\tens1 : M\tens X\tens N\tens Y \to M\tens N\tens X\tens Y
\end{equation*}
for $M,N\in {}^H\cv$, $X,Y\in\cv$, gives a monoidal equivalence
\begin{equation*}\label{e166.2b}
(\ot, \psi, r_I^{-1}) : ({}^H\cv\tens\cv, \unten, I\odot I) \to
({}^{\bar H}\cv, \tens, I) .
\end{equation*}
\end{proposition}

Each of the four isomorphisms
\begin{multline*}
c_{\pm\pm}: (M\odot X)\unten(N\odot Y) = (M\tens N)\odot(X\tens Y) \\
\xra{R^{\pm1}\odot c^{\pm1}} (N\tens M)\odot(Y\tens X) =
(N\odot Y)\unten(M\odot X)
\end{multline*}
gives a braiding in $({}^H\cv\tens\cv, \unten, I\odot I)$. Thus, the
category of $\bar H$-comodules is braided as well.

\begin{corollary}
The category $({}^{\bar H}\cv, \tens, I)$ has braidings $c_{\pm\pm}'$ such
that
\begin{multline*}
c_{\pm\pm}'= \bigl( (M\tens X)\tens (N\tens Y) \xra{1\tens c\tens1}
M\tens N\tens X\tens Y \\
\xra{R^{\pm1}\tens c^{\pm1}} N\tens M\tens Y\tens X
\xra{1\tens c^{-1}\tens1}  (N\tens Y)\tens (M\tens X) \bigr)
\end{multline*}
for $M,N\in {}^H\cv$, $X,Y\in\cv$.
\end{corollary}

\section{Ribbon Hopf coalgebras}
Let us assume that $\cv$ is braided. Fix the isomorphism $\zeta = u_1^2 : X
\to X^{\vee\vee}$. We want to discuss ribbon structures in the categories
${}^H\cv$. A ribbon structure in $\cv$ is not required.

\subsection{Ribbon twists}
\begin{definition}
A \emph{ribbon Hopf coalgebra} is a quasitriangular Hopf coalgebra $H$
equipped with a linear form $\Theta: \bar H \to I \in \vhat$ such that
\begin{equation}\label{e171a}
\bigl( I \xra{\bar\eta} \bar H \xra\Theta I \bigr) = \id_I ,
\end{equation}
\begin{equation}\label{e171b}
\bigl( H^{12} \xra{c^{12}} H^{21} \xra{\gamma_\zeta^{12}} H^{12} \xra\Theta
I \bigr) = \Theta ,
\end{equation}
\begin{equation}\label{e171c}
\begin{split}
&\bigl( \bar H\tens\bar H \xra{\bar\Delta\tens\bar\Delta} \bar H\tens\bar
H\tens\bar H\tens\bar H \xra{\bar\Delta\tens\Theta\tens\bar H\tens\Theta}
\bar H\tens\bar H\tens I\tens\bar H\tens I \\
&\hspace*{3cm} \xra\sim \bar H\tens\bar H\tens\bar H \xra{\bar H\tens c}
\bar H\tens\bar H\tens\bar H \xra{\bar H\tens\bar\Delta\tens\bar H} \\
&\hspace*{3cm} \bar H\tens\bar H\tens\bar H\tens\bar H
\xra{R_+\tens R_-} I\tens I \xra\sim I \bigr) \\
&= \bigl( \bar H\tens\bar H \xra{\bar m} \bar H \xra\Theta I \bigr) ,
\end{split}
\end{equation}
\begin{multline}\label{e171d}
\bigl( H_{12} \xra\sim H_{1'2}\tens I_{1''} \xra{\Delta_{1'1''2}}
H_{1'1''}\tens H_{1'''2} \xra{\Theta\tens H} I_{1'}\tens H_{1''2} \xra\sim
H_{12} \bigr) \\
= \bigl( H_{12} \xra\sim H_{12''}\tens I_{2'} \xra{\Delta_{12'2''}}
H_{12'}\tens H_{2''2'''} \xra{H\tens\Theta} H_{12'}\tens I_{2''} \xra\sim
H_{12} \bigr) .
\end{multline}
\end{definition}

\begin{theorem}
\textup{(a)}
The category ${}^H\cv$ of comodules over a quasitriangular Hopf coalgebra
$H$ is ribbon if and only if $H$ is ribbon.

\textup{(b)}
A linear form $\Theta: \bar H \to I\in \vhat$ satisfying the equations
\eqref{e171a}--\eqref{e171d} determines a ribbon twist in ${}^H\cv$
by the formula
\[ \theta_X = \bigl( X \xra{\bar\delta} \bar H\tens X \xra{\Theta\tens X}
I\tens X \simeq X \bigr) \in \cv. \]

\textup{(c)}
A ribbon twist $\theta$ in ${}^H\cv$ determines the unique functional
$\Theta$ with the properties \eqref{e171a}--\eqref{e171d} via
the diagram
\[
\begin{CD}
X\tens X\pti @>\theta\tens X\pti>> X\tens X\pti \\
@V\ot\ii_XVV  @VV\ev V \\
\bar H @>\exists\Theta>\hphantom{\theta\tens X\pti}> I
\end{CD}
\qquad .
\]
\end{theorem}

\subsubsection*{Acknowledgements} I would like to thank Yu.~Bespalov
and A.~Sudbery for fruitful discussions and attention to this work.

\end{document}